\documentclass[10pt,twoside,reqno]{amsart}
 \usepackage[top=1in, bottom=1in, left=1in, right=1in]{geometry}
\usepackage{graphicx}
\usepackage{amsmath,latexsym,amssymb}
\usepackage{amsfonts}
\usepackage{mathtools,esint}
\usepackage{amsmath}
\usepackage{amscd,amsthm}
\usepackage{stmaryrd}
\usepackage{fancyhdr}
\usepackage{commath,enumerate}
\usepackage{indentfirst, hyperref}
\usepackage{caption, subcaption}
\usepackage{bbold}
\usepackage{pgf,tikz}
\usepackage{mathrsfs}
\usetikzlibrary{arrows,calc,positioning}
\usepackage{float}
\numberwithin{figure}{section}
\usepackage[bottom]{footmisc}

\numberwithin{equation}{section}
\newtheorem{theorem}{Theorem}[section]
\newtheorem{lemma}[theorem]{Lemma}
\newtheorem{proposition}[theorem]{Proposition}

\theoremstyle{definition}

\newcommand{\Rm}[1]{
  \textup{\uppercase\expandafter{\romannumeral#1}}
}

\newcommand{\F}{\mathcal{F}}
\newcommand{\C}{\mathcal{C}}

\newcommand{\gl}{{\gamma_l}}
\newcommand{\gh}{{\gamma_h}}

\newcommand{\Nc}{\mathcal{N}}

\renewcommand{\P}{\mathcal{P}}
\newcommand{\R}{\mathbb{R}}

\newcommand{\Tb}{\mathbf{T}}

\renewcommand{\S}{\mathcal{S}}
\newcommand{\Z}{\mathbb{Z}}
\newcommand{\etab}{\boldsymbol{\eta}}
\newcommand{\cutoffxi}{\mathfrak{d}}

\newcommand{\cutoff}{\mathfrak{b}}

\def\vp{\varphi}
\def\ve{\varepsilon}
\def\px{\partial_x}
\def\p{\partial}

\newcommand{\diff}{\,\mathrm{d}}

\usepackage{mathtools}

\author{Fangchi Yan}
\address{Department of Mathematics, Virginia Tech}
\email{fyan1@alumni.nd.edu}
\author{Qingtian Zhang}
\address{Department of Mathematics, West Virginia University}
\email{qingtian.zhang@mail.wvu.edu}
\title[quasilinear mkdv]{Global solutions of quasi-linear Hamiltonian mKdV equation}
\date{\today}

\begin{document}

\begin{abstract}

We study the initial value problem of quasi-linear Hamiltonian mKdV equations. Our goal is to prove the global-in-time existence of a solution given sufficiently smooth, localized, and small initial data. To achieve this, we utilize the bootstrap argument, Sobolev energy estimates, and the dispersive estimate. This proof relies on the space-time resonance method, as well as a bilinear estimate developed by Germain, Pusateri, and Rousset. 


\end{abstract}

\maketitle

\tableofcontents

\section{Introduction}
This paper focuses on the Cauchy problem of the quasilinear modified Korteweg–de Vries (mKdV) equation, given by:
\begin{subequations}
\label{qmkdv}
\begin{align}
\label{qmkdv:eqn}
&\vp_t+\px^{3}\vp
+
\px(\vp^
{3}
)
+
\px(c(\vp)\px(c(\vp)\px \vp))
=
0,
\quad
x\in\R,\\
\label{qmkdv:ic}
&\vp(0, x)
=
\vp_0(x),
\quad
x\in\R.
\end{align}
\end{subequations}
{
Here, $c(\vp)$ is a function of $\vp$, and this equation belongs to a broader category of quasilinear KdV-type equations, which can be expressed as:
\begin{equation}
\label{quas-kdv-eqn}
\vp_t+\px^{3}\vp
+
\px(\vp^{k})
+
\px(c(\vp)\px(c(\vp)\px \vp))
=
0.
\end{equation}
 Equation \eqref{quas-kdv-eqn} possesses a Hamiltonian structure, which can be written as
}
\[
\vp_t=\left[\frac{\delta\mathcal{H}(\vp)}{\delta\vp}\right]_x,
\]
where the Hamiltonian functional is given by 
$\mathcal H(\vp)=-\frac1{k+1}\vp^{k+1}+\frac12(c^2(\vp)+1)\vp_x^2$.

When $c=0$, the equation \eqref{quas-kdv-eqn} can be referred to as the KdV equation when $k=2$, the modified KdV (mKdV) equation when $k=3$, or the generalized KdV (gKdV) equation for any positive integer $k\ge 4$. 

There is a long history of the study of  KdV type equations, and we will mention some related results here. Bona and Smith \cite{BS75} showed that local well-posedness in $H^s(\mathbb{R})$, $s>\frac{3}{2}$, and global well-posedness in $H^s(\mathbb{R})$, $s\geq 2$, hold for the KdV equation. Kato \cite{Kato1983} lowered the requirement of $s$ by proving that if the initial data is in $H^0\cap L^2(e^{2bx} dx)$, the Cauchy problem of KdV equation is globally well-posed.  
Additionally, Kenig, Ponce, and Vega \cite{kpv1991} showed that the KdV equation is locally well-posed in $H^s(\mathbb{R})$, where $s>\frac34$.
Bourgain \cite{b1993-kdv} later proved the well-posedness of the KdV equation in $H^s(\mathbb{T})$ and $H^s(\mathbb{R})$ for $s\geq 0$ using bilinear estimates in the spaces $X^{s,b}$ introduced there. Utilizing these spaces, Kenig, Ponce, and Vega \cite{kpv96} proved the local well-posedness of KdV  on the line for $s>-3/4$ and on the circle for $s> -1/2$. In addition, Colliander, Keel, Staffilani, Takaoka, and Tao \cite{CKSTT03} proved the global well-posedness of the KdV equation for the same range of Sobolev exponents.
Moreover, Guo \cite{Guo2009} and Kashimoto \cite{Kis2009} independently proved the  global well-posedness of the KdV equation on the line for the critical Sobolev exponent $s=-3/4$. 
As for the mKdV equation, Tao \cite{tao2001} proved its well-posedness for $s\geq 1/4$ using trilinear estimates in Bourgain spaces. An earlier proof, which is part of gKdV well-posedness, can be found in \cite{kpv93b}.

The main technique for proving the global well-posedness of KdV type equations is the ``I-method," which was developed by Tao and his collaborators \cite{tao2001,CKSTT03,CKSTT04}. The ``I-method" involves decomposing the equation into interaction and non-interaction parts and then applying frequency modulation and Littlewood-Paley decomposition techniques to estimate the interaction terms. This iterative process of decomposition and estimation allows one to get multilinear estimates in Bourgain-type spaces and
 derive global well-posedness in the Bourgain-type spaces. 

However, the Bourgain-type spaces mentioned above are a strict subset of the Hadamard spaces (natural solution spaces) given by
\begin{equation}
\label{Hadamard-space}
X = C(\mathbb{R}; H^s(\mathbb{R})).
\end{equation}
The well-posedness results mentioned earlier can only guarantee uniqueness of solutions in a specific subset of Hadamard spaces $X$, for solutions with lower regularity. Such well-posedness is referred to as conditional well-posedness. On the other hand, unconditional well-posedness is when uniqueness of solutions is guaranteed in the entire Hadamard space.

Over the last 25 years, unconditional well-posedness (WP) for nonlinear dispersive equations has been an intriguing topic. Kato \cite{Kt1995} first established unconditional WP for the nonlinear Schr\"odinger equation (NLS), and Zhou \cite{Zy1997} established it for the KdV equation. For the mKdV equation, unconditional WP was achieved through the normal form approach by Kwon and Oh \cite{KO12}. For further research on the topic, we recommend \cite{GKO2013,KOY20,GL2021,Kis2022-1,Kis2022-2} to the reader.

Another key feature of KdV and mKdV is their integrability. The modern theory of integrable systems began in 1967 by Gardner, Greene, Kruskal and Miura, who introduced the inverse scattering transform to solve KdV equation \cite{GGKM}. In 1968, Lax \cite{Lax} found a new method to solve KdV equation based on two operators, which was called Lax pairs in later references. Afterwards Lax pairs were found in many other equations or systems. See the monograph by Ablowitz and Clarkson \cite{AbCl91} for the detailed discussion. The integrability of KdV and mKdV is also a great benefit in the study of the long-time behavior. 
Segur and Ablowitz provided a complete description of the leading asymptotics for the solution of the Cauchy problem for KdV, MKdV, and Sine-Gordon equations in their papers \cite{SA77, SA81}.
Similarly, Its introduced a method using the Riemann-Hilbert problem to obtain complete asymptotic behavior for solutions to the NLS equation in 1981, as described in \cite{Its81}.
Deift and Zhou then built on this work by developing the nonlinear steepest descent method for analyzing the asymptotics of oscillatory Riemann-Hilbert problems.

Recently, there are several results studying the long-time behavior of the solutions to mKdV equation without using integrable structure.  Hayashi and Naumkin \cite{HaNa99} using the method of factorization of operators proved the global existence and asymptotic behavior for small localized 
solutions. Germain-Pusateri-Rousset \cite{GPR16} used the space-time resonance method proved the asymptotic stability of zero solution and the solitons for mKdV equation.  Harrop-Griffiths \cite{HaGr16} obtained the similar result for zero solution by using wave-packets analysis. Chen and Liu \cite{CL22} using the nonlinear steepest descent method of Deift and Zhou proved the long-time asymptotics in weighted Sobolev spaces. In \cite{CL21}, Chen and Liu studied the soliton resolution for focusing mKdV equation.

The quasilinear KdV type equation is very different from its semi-linear counterparts. First, due to the quasilinear structure, the technique based on the Duhamel formula, which is a convenient way to express solutions of semi-linear KdV type equations, will not apply. Secondly, only very few cases are known to be integrable \cite{Mar10, Mar10b}. Due to its quasilinear structure, the method of integrable system cannot be easily adapted. Based on these two reasons, it seems interesting to achieve a better understanding of the solutions to quasilinear KdV type equations.

Currently, the results about quasilinear KdV type equations are fewer compared with the above semi-linear counter-part. The local-in-time solutions of Cauchy problems of quasilinear KdV equations with more general nonlinearities are proved in \cite{Ak2013, CKS92, AAW2019, LPS17}, under the nondegenerate assumption of the third-order coefficient. In the general setting, one assumption is required to control the strength of backwards diffusion. When the degeneracy happens, \cite{GHM2019} proved the local existence and uniqueness of solutions. In our equation \eqref{qmkdv}, there is no backwards diffusion or degeneracy. Its local well-posedness is proved by Mietka \cite{Mie17} for $s\geq 4$ and improved later by Iandoli \cite{Ian22} for $s>3.5$.
 
This paper is one step further to understanding the quasilinear KdV-type equations. 
To state our main result, we introduce the following notations:
\begin{equation}
\label{h-def}
h(t, x) = e^{t\px^3}\vp(t, x), \qquad \hat{h}(t, \xi) = e^{-it\xi^3}\hat{\vp}(t, \xi),
\end{equation}
where $h$ is obtained by removing the action of the linearized evolution group on $\vp$. We also define the scaling vector field as
\begin{equation}
\label{scaling-operator}
\S = x\p_x+3t\p_t.
\end{equation}
In this paper, we define $\delta$ to be a small positive number, specifically
\begin{equation}
\label{delta-choice}
\delta
\doteq
10^{-3}.
\end{equation} 
%
%
For $s\ge 12$ and $\delta$ given by \eqref{delta-choice}, define $p_0$:
\begin{equation}
\label{p0-choice}
p_0
\doteq
\frac1{10}\delta
=
10^{-4}.
\end{equation}

We define the energy of the solution $\vp$ by
\begin{align}
\label{Def_energy}
E[\vp]=\|\px^{-1}\vp(t)\|_{L^2}^2+\|\vp(t)\|_{H^s}^2+\|\px^{-1}\S\vp(t)\|_{L^2}^2+\|\S\vp(t)\|_{L^2}^2+\|\xi\partial_\xi\hat h(t)\|_{L^2_\xi}^2+\|\partial_\xi\hat h(t)\|_{L^2_\xi}^2.
\end{align}
Compared with \cite{GPR16}, we included the low-frequency bound $\|\px^{-1}\vp(t)\|_{L^2}$ based on a symmetry consideration. In gain, we get a better estimate on $\|\px^{-1}\S\vp(t)\|_{L^2}$.

Finally, 
{
letting $\gl=0$ and $\gh=\frac52$,
}
 we define the $Z$-norm as follows:
\begin{equation}\label{Def_Z}
\|\vp\|_{Z}=\big\|(|\xi|^\gl+|\xi|^\gh) \hat \vp(\xi)\big\|_{L^\infty_\xi},
\end{equation}
%
Using the above notations, our main result can be stated as follows:
\begin{theorem}
\label{Thm:glo}
Assume
{$s\ge 12$ and}
 $c: \R\to \R$ is in $C^{s+2}$ satisfying $c(0)=0$.  
There exists a constant $0\le \varepsilon \ll 1$ such that for any $0\le\varepsilon_0\le \varepsilon$,
$\vp_0\in H^s(\R)$ satisfying
\begin{align}
\label{initial-assump}
E[\vp_0]+\|\vp_0\|_{Z}^2
\le
\varepsilon_0^2,
\end{align}
 there exists a unique global solution $\vp \in C([0,\infty); H^s(\R))$ of \eqref{qmkdv}. Moreover, this solution satisfies 
\begin{equation}
(1+t)^{-2p_0}E[\vp](t)+\|\vp(t)\|_{Z}^2 \lesssim \varepsilon_0^2,
\end{equation}
where 
{
$p_0=10^{-4}$ is the small number defined by \eqref{p0-choice}.
}
\end{theorem}
%
%

The proof of Theorem \ref{Thm:glo} relies on the space-time resonance method developed by Germain, Masmoudi, and Shatah \cite{Germain,GMS09, GMS12}, in conjunction with $Z$-norm estimates, which are weighted $L^\infty_\xi$-norms developed by Ionescu and collaborators \cite{CGI19, DIP17, DIPP17, IP12, IP13, IP14, IP15, IPu16}. This method has also been employed to solve global solutions of vortex front problems \cite{HSZ21, HSZ20p, YZ22p}. The core idea of this approach is to prove the decay of $L^\infty_x$ norm by estimating the weighted $L^\infty_\xi$-norms. The modified scattering and resonant analysis will be employed to establish these weighted $L^\infty_\xi$-norms. For the KdV equation, the dispersion relation $p(\xi)=\xi^3$ has $p''(0)=0$, indicating that $0$ frequency is a stationary phase. The linear dispersive estimate alone may only provide $t^{-1/3}$ decay rate, which may not be sufficient for the cubic nonlinearity. 
To overcome this challenge, 
we use the bilinear estimate \eqref{dec1} by 
Germain-Pusateri-Rousset  \cite{GPR16}. 

{\it Structure of the paper.} 
Section \ref{sec:prelim} collects some useful lemmas and notations. In Section \ref{Sec:structure}, we analyze the structure of the equation.  The estimates of $\|\vp(t)\|_{H^s}^2$, $\|\px^{-1}\vp(t)\|_{L^2}^2$ and $\|\S\vp(t)\|_{L^2}$, $\|x h(t)\|_{L^2}$ are in Section \ref{priori-energy-est}. 
The dispersive estimate and the sharp point-wise decay are proved in Section \ref{sec:dis}.  
  Section \ref{sec:Znorm} is for Z-norm estimate. 

\section{Preliminaries and Strategy of the proof}
\label{sec:prelim}
\subsection{Notations and useful tools}
We denote the Fourier transform of $f \colon \R\to \C$ by $\hat f \colon \R\to \C$, where $\hat f= \F f$ is given by
\[
f(x)=\int_{\R} \hat f(\xi) e^{i\xi x} \diff\xi,  \qquad \hat f(\xi)=\frac1{2\pi} \int_{\R}f(x) e^{-i\xi x}\diff{x}.
\]
For $s\in \R$, we denote by $H^s(\R)$ the space of Schwartz distributions $f$ with $\|f\|_{H^s} < \infty$, where
\[
\|f\|_{H^s} = \left[\int_\R \left(1+|\xi|^2\right)^s |\hat{f}(\xi)|^2\, \diff{\xi}\right]^{1/2}.
\]
Throughout this paper, we use $A\lesssim B$ to mean there is a constant $C$ such that $A\leq C B$, and $A\gtrsim B$ to mean there is a constant $C$ such that $A\geq C B$. We use $A\approx B$ to mean that $A\lesssim B$ and $B\lesssim A$.

Let $\psi \colon\R\to [0,1]$ be a smooth function supported in $[-8/5, 8/5]$ and equal to $1$ in $[-5/4, 5/4]$.
For any $k\in \mathbb Z$, we define
\begin{align}
\label{defpsik}
\begin{split}
\psi_k(\xi)&=\psi(\xi/2^k)-\psi(\xi/2^{k-1}), \qquad \psi_{\leq k}(\xi)=\psi(\xi/2^k),\qquad \psi_{\geq k}(\xi)=1-\psi(\xi/2^{k-1}),\\
\tilde\psi_k(\xi)&=\psi_{k-1}(\xi)+\psi_k(\xi)+\psi_{k+1}(\xi).
\end{split}
\end{align}
So we have the homogeneous dyadic decomposition
\[
1=\sum\limits_{k\in\mathbb Z}\psi_k(\xi),
\]
and the non-homogeneous dyadic decomposition
\[
1=\sum\limits_{k\in\mathbb N}\psi_k(\xi)+\psi_{\leq 0}(\xi).
\]
We denote the sub-index of the non-homogeneous dyadic decomposition as $\{\leq 0\}\cup \mathbb N$.

We denote by $P_k$, $P_{\leq k}$, $P_{\geq k}$, and $\tilde{P}_k$  the Fourier multiplier operators with symbols $\psi_k, \psi_{\leq k}, \psi_{\geq k}$, and $\tilde{\psi}_k$, respectively. Notice that $\psi_k(\xi)=\psi_0(\xi/2^k)$, $\tilde\psi_k(\xi)=\tilde\psi_0(\xi/2^k)$. 
When there is no ambiguity, we denote $P_k\vp$ in short as $\vp_k$.

It is easy to check that
\begin{equation}
\label{psi-L2}
 \|\psi_k\|_{L^2}\approx  2^{k/2}, \qquad  \|\psi_k'\|_{L^2}\approx 2^{-k/2}.
\end{equation}

We will need the following interpolation lemma, whose proof can be found in \cite{IPu16}.
\begin{lemma}\label{interpolation}
For any $k\in\mathbb N$ and $f\in L^2(\R)$, we have
\[
\|\widehat{P_kf}\|_{L^\infty}^2\lesssim \|P_k f\|_{L^1}^2\lesssim 2^{-k}\|\hat f\|_{L^2_\xi}\left[2^k\|\partial_\xi\hat f\|_{L^2_\xi}+\|\hat f\|_{L^2_\xi}\right].
\]
\end{lemma}

Define the symbol class
\[
S^\infty:=\{m: \R^d\to\C: m \text{ is continuous and }  \|m\|_{S^\infty}:=\|\F^{-1}(m)\|_{L^1}<\infty\}.
\]

\begin{lemma}[Algebraic property of $S^\infty$]
If $m$, $m'\in S^\infty$, then $m\cdot m'\in S^\infty$.
\end{lemma}

Given $m \in S^\infty(\R^{n})$, we define a multilinear operator $M_m$ acting on Schwartz functions  $f_1,\dots, f_m \in \mathcal{S}(\R)$ by
\[
M_{m}(f_1,\dotsc,f_n)(x)=\int_{\R^{n}} e^{ix(\xi_1+\dotsb+\xi_n)}m(\xi_1, \dotsc, \xi_n)\hat f_1(\xi_1) \dotsm \hat f_n(\xi_n)\diff{\xi_1} \dotsm \diff{\xi_n}.
\]

\begin{lemma}[Estimates of multilinear Fourier integral operators]\label{multilinear}
(i)~ Suppose that $1<p_1, \dotsc, p_n\leq \infty$, $0<p<\infty$,  satisfy
\[
\frac1{p_1}+\frac1{p_2}+ \dotsb +\frac1{p_n}=\frac1p,
\]
for every nonempty subset $J\subset \{1,2, \dotsc, n\}$. If $m \in S^\infty(\R^{n})$, then
\[
\|M_{m}\|_{L^{p_1}\times \dotsb \times L^{p_n}\to L^p}\lesssim \|m\|_{S^\infty}.
\]
(ii)~ Assume $p,q,r\in[1,\infty]$ satisfy $\frac1p+\frac1q=\frac1r$, and $m\in S^\infty$. Then for any $f,g\in L^2(\R)$, 
\[
\left|\int_{\R^2}m(\xi,\eta)\hat f(\xi)\hat g(\eta)\hat h(-\xi-\eta)\diff \xi\diff\eta\right|\lesssim \|m\|_{S^\infty}\|f\|_{L^p}\|g\|_{L^q}\|h\|_{L^r}.
\]
\end{lemma}

\begin{lemma}[Gagliardo-Nirenberg inequality]\label{GN-inter-ine}
If $u\in L^q(\R)$ and $\px^m u\in L^{r}(\R)$, $1\le q, r\le \infty$ and $0\le j\le m$,
then we have $\px^{j} u\in L^p(\R)$, where  $p$ satisfies
\begin{align}
\label{jp-cond}
\frac{1}{p}
=
j
+
\left(
\frac{1}{r}
-
m
\right)
\alpha
+
\frac{1-\alpha}{q},
\quad
\frac{j}{m}
\le
\alpha
\le
1.
\end{align}
More precisely, we have the following inequality
\begin{equation}
\label{GN-ine}
\|\px^j u\|_{L^p}
\le
C
\|\px^m u\|_{L^r}^{\alpha}
\cdot
\|u\|_{L^q}^{1-\alpha},
\end{equation}
where $C=C(m,n,j,q,r,\alpha)$.
\end{lemma}

Denote $|f(t, x)|_s=| f(t,x)|+|\px f(t,x)|+\cdots+|\px^sf(t,x)|$.

We define the $B^{a,b}$ semi-norm as 
\begin{equation}\label{bnorm}
\|f\|_{B^{a,b}}=\sum\limits_{j\in\Z}(2^{aj}+2^{bj})\|P_jf\|_{L^\infty}.
\end{equation}
Here $a\leq b$ and $a,b$ are the indices of higher and lower frequencies. Compared with the $L^\infty$ norm, $B^{a,b}$ satisfies one property not shared with the $L^\infty$ norm, that is
\[
\|f\|_{B^{a,b}}\approx \sum\limits_{j\in\Z}\|P_jf\|_{B^{a,b}}.
\]
When $a=b$, we write it in short as $B^a:=B^{a,a}.$
Notice that $\||\px|^af\|_{W^{b-a,\infty}}\lesssim \|f\|_{B^{a,b}}$.

\subsection{Strategy of the proof}
The theorem \ref{Thm:glo} follows from the global energy estimate by the following boot-strap argument.
\begin{proposition}
\label{Bootstrap-prop}
Let $T>1$ and suppose that $\vp \in C([0,T]; H^s)$ is a solution of \eqref{qmkdv}, where the initial data satisfies \eqref{initial-assump} for some $0\le\varepsilon_0\le 1$. 
If there exists $\varepsilon_0\ll\varepsilon_1\lesssim \varepsilon_0^{2/3}$ such that the solution satisfies
\begin{equation}
\label{bootstrap-assump}
(1+t)^{-2p_0}E[\vp](t)+\|\vp(t)\|_{Z}^2\le
\varepsilon_1^2, \quad \text{for every }  t\in[0,T],
\end{equation}
where $E[\vp]$ and Z-norm are defined in \eqref{Def_energy} and \eqref{Def_Z} separately,
 then the solution satisfies an improved bound
\begin{equation}
(1+t)^{-2p_0}E[\vp](t)+\|\vp(t)\|_{Z}^2 \le
\varepsilon_0^2 \quad \text{for every }  t\in[0,T].
\end{equation}
\end{proposition}

\begin{lemma}[Energy estimate]\label{energyLem}
Recall the definition of $E[\vp](t)$ in \eqref{Def_energy}. 
Under the bootstrap assumptions, the solution of \eqref{qmkdv} satisfies the energy inequality
\begin{equation}\label{EnergyIneq}
E[\vp](t)\lesssim E[\vp_0]+ \int_0^t F(\|\vp\|_{W^{2,\infty}})\||\vp|_2|\px\vp|_2\|_{L^\infty}E[\vp](\tau)d\tau.
\end{equation}
\end{lemma}
This Lemma directly leads to the improved bound of $E[\vp](t)$ in Proposition \ref{Bootstrap-prop}. The proof of Lemma \ref{energyLem} is in Section \ref{priori-energy-est}.

\begin{lemma}
[Sharp pointwise decay]
\label{lem:decay}
Under the bootstrap assumptions, 
\begin{align}\label{dec2}
\||\px|^{1/2}\vp(t,x)\|_{L^\infty}+\|\vp_x(t,x)\|_{L^\infty}+\|\vp_{xx}(t,x)\|_{L^\infty}+\|\vp_{xxx}(t,x)\|_{L^\infty}\lesssim t^{-1/2}\ve_1,
\end{align}
\begin{align}
\label{dec1}
\||\vp|_2|\px\vp|_2\|_{L^\infty}\lesssim\varepsilon_1^2 (1+t)^{-1}.
\end{align}
\end{lemma}

\section{Structure of the equation}\label{Sec:structure}
From \eqref{qmkdv:eqn}, we  rewrite the quasilinear mKdV equation as follows:

\begin{equation}
\label{qmkdv:Eq}
\vp_t
+
\p_x^3\vp
=
-\Nc(\vp),
\end{equation}
where
\begin{align*}
\Nc(\vp)
=
\px(\vp^3)
+
\px(c(\vp)\px(c(\vp)\px \vp)).
\end{align*}

Since $c(0)=0$, we write $c$ as $c(\vp)=c'(0)\vp+\frac12c''(0)\vp^2+c_3(\vp)$, where $c_3$ satisfies
\[
c_3(0)=c'_3(0)=c''_3(0)=0.
\]
Then we have
\begin{align*}
\px(c(\vp)\px(c(\vp)\px \vp))=&\px\left([c'(0)\vp+\frac12c''(0)\vp^2+c_3(\vp)]\px([c'(0)\vp+\frac12c''(0)\vp^2+c_3(\vp)]\px \vp)\right)
\end{align*}
we  rewrite the nonlinearity $\Nc(\vp)$ as
\begin{align*}
\Nc(\vp)
=
\Nc_3(\vp)
+\Nc_4(\vp)+
\Nc_{\ge 5}(\vp),
\end{align*}
where
the cubic nonlinearity is given by
\begin{align*}
\Nc_3(\vp)
=
\px (\vp^3)
+
\px
\left[
\alpha_2\vp^2\vp_{xx}
+
\alpha_2\vp(\vp_{x})^2
\right],
\quad
\text{where}
\,\,
\alpha_2
\doteq
[c'(0)]^2,
\end{align*}
and the higher order nonlinearity is given by
\begin{align*}
\Nc_{ 4}(\vp)
=\alpha_3
\px
\left[
\vp^2\px(\vp\px \vp)+\vp\px(\vp^2\px \vp)
\right],\qquad \text{where}\,\, \alpha_3=\frac12c''(0)c'(0),
\end{align*}
\begin{align*}
\Nc_{\geq5}(\vp)
=&\px\left(c_3(\vp)\px([c'(0)\vp+\frac12c''(0)\vp^2+c_3(\vp)]\px \vp)\right)+
\px\left([c'(0)\vp+\frac12c''(0)\vp^2+c_3(\vp)]\px(c_3(\vp)\px \vp)\right)\\
&+\px\left(\frac12c''(0)\vp^2\px(\frac12c''(0)\vp^2\px \vp)\right).
\end{align*}
For a shorter notation, we denote $\Nc_{\geq 4}=\Nc_{4}+\Nc_{\geq 5}$.

Taking Fourier transform, we get
\begin{equation}
\label{vp-FT-eqn}
\p_t\hat{\vp}(t, \xi)
-
i\xi^3
\hat{\vp}(t, \xi)
=
-
\widehat{\Nc(\vp)}(t, \xi)
=
-
\widehat{\Nc_3(\vp)}(t, \xi)
-
\widehat{\Nc_{\ge 4}(\vp)}(t, \xi)
\end{equation}
Next, we write $\widehat{\Nc_3(\vp)}(t, \xi)$ explicitly. In fact, using $\widehat{\p_x \vp}(\xi)=i\xi\hat{\vp}(\xi)$ we get
$$
-\widehat{\Nc_3(\vp)}(t, \xi)
=
i\xi
\iint
(\alpha_2\eta_1^2+\alpha_2\eta_1\eta_2-1)
\hat{\vp}(\xi-\eta_1-\eta_2,t)
\hat{\vp}(\eta_1,t)
\hat{\vp}(\eta_2,t)
d\eta_1
d\eta_2.
$$
In addition, letting $\eta_3=\xi-\eta_1-\eta_2$, we get 
\begin{align}
\label{N-phi-3}
-\widehat{\Nc_3(\vp)}(t, \xi)
=
i\xi
\iint
(\alpha_2\eta_1^2+\alpha_2\eta_1\eta_2-1)
\hat{\vp}(\eta_3,t)
\hat{\vp}(\eta_1,t)
\hat{\vp}(\eta_2,t)
d\eta_1
d\eta_2.
\end{align}
Thus, by using symmetry we define
\begin{align}
\label{T1-def-sym}
{\bf T_1}(\eta_1,\eta_2,\eta_3)
\doteq
&\frac16
\Big[
(\alpha_2\eta_1^2+\alpha_2\eta_1\eta_2-1)
+
(\alpha_2\eta_1^2+\alpha_2\eta_1\eta_3-1)
+
(\alpha_2\eta_2^2+\alpha_2\eta_2\eta_1-1)
\\
+&
(\alpha_2\eta_2^2+\alpha_2\eta_2\eta_3-1)
+
(\alpha_2\eta_3^2+\alpha_2\eta_3\eta_1-1)
+
(\alpha_2\eta_3^2+\alpha_2\eta_3\eta_2-1)
\Big]
\notag
\\
=&
\frac{\alpha_2}3(\eta_1^2+\eta_2^2+\eta_3^2+\eta_1\eta_2+\eta_1\eta_3+\eta_2\eta_3)-1.
\notag
\end{align}
Furthermore,  by letting $\eta_3=\xi-\eta_1-\eta_2$ we have 
\begin{align}
\label{T1-def}
{\bf T_1}(\eta_1,\eta_2,\xi-\eta_1-\eta_2)
=
\frac{\alpha_2}3(\eta_1^2+\eta_2^2+\xi^2+\eta_1\eta_2-\eta_1\xi-\eta_2\xi)-1.
\end{align}
And we can rewrite \eqref{N-phi-3} as
\begin{align}
\label{N-phi-3-eqn2}
-\widehat{\Nc_3(\vp)}(t, \xi)
=&
i\xi
\iint
{\bf T_1}(\eta_1,\eta_2,\eta_3)
\hat{\vp}(\eta_1,t)
\hat{\vp}(\eta_2,t)
\hat{\vp}(\eta_3,t)
d\eta_1
d\eta_2
\\
=&
i\xi
\iint
{\bf T_1}(\eta_1,\eta_2,\xi-\eta_1-\eta_2)
\hat{\vp}(\eta_1,t)
\hat{\vp}(\eta_2,t)
\hat{\vp}(\xi-\eta_1-\eta_2,t)
d\eta_1
d\eta_2.
\nonumber
\end{align}
Therefore, substituting \eqref{N-phi-3-eqn2} into \eqref{vp-FT-eqn}, we get 
\begin{equation}
\label{vp-FT-fin}
\p_t\hat{\vp}(t, \xi)
-
i\xi^3
\hat{\vp}(t, \xi)
=
i\xi
\iint
{\bf T_1}(\eta_1,\eta_2,\xi-\eta_1-\eta_2)
\hat{\vp}(\eta_3,t)
\hat{\vp}(\eta_1,t)
\hat{\vp}(\eta_2,t)
d\eta_1
d\eta_2
-
\widehat{\Nc_{\ge 4}(\vp)}(t, \xi).
\end{equation}
%
Using definition \eqref{h-def}, equation \eqref{vp-FT-eqn} becomes
\begin{align}
\label{f-eqn}
\p_t \hat{h}(t, \xi)
=
-e^{-i\xi^3 t}
\Big[\p_t\hat{\vp}(t, \xi)
-
i\xi^3
\hat{\vp}(t, \xi)\Big]
\end{align}
Now, using \eqref{vp-FT-fin}, we get
\begin{align}
\label{h-FT-fin}
\p_t \hat{h}(t, \xi)
=&
\iint
e^{-it\Phi}
i\xi
\,
{\bf T_1}(\eta_1,\eta_2,\eta_3)
\hat{h}(\eta_3,t)
\hat{h}(\eta_1,t)
\hat{h}(\eta_2,t)
d\eta_1
d\eta_2
-e^{-i\xi^3 t}\widehat{\Nc_{\ge 4}(\vp)}(t, \xi)
\\
=&
\iint
e^{-it\Phi}
i\xi
\,
{\bf T_1}(\eta_1,\eta_2,\xi-\eta_1-\eta_2)
\hat{h}(\xi-\eta_1-\eta_2,t)
\hat{h}(\eta_1,t)
\hat{h}(\eta_2,t)
d\eta_1
d\eta_2
-e^{-i\xi^3 t}\widehat{\Nc_{\ge 4}(\vp)}(t, \xi),
\nonumber
\end{align}
where $\Tb_1$ is given by \eqref{T1-def} and 
with
\begin{equation}
\label{defPhi}
\Phi
=
\Phi(\eta_1,\eta_2,\xi-\eta_1-\eta_2)
\doteq
\xi^3-(\xi-\eta_1-\eta_2)^3
-
\eta_1^3
-
\eta_2^3
=
3(\eta_1+\eta_2)(\xi-\eta_1)(\xi-\eta_2).
\end{equation}
Furthermore, taking derivative we get
\begin{align*}
\p_{\eta_1}\Phi
=
3(\xi-\eta_2)(\xi-2\eta_1-\eta_2),
\quad
\p_{\eta_2}\Phi
=
3(\xi-\eta_1)(\xi-2\eta_2-\eta_1).
\end{align*}
Hence, we obtain
$$
\p_{\eta_1}\Phi
=
\p_{\eta_2}\Phi
=
0
\Longleftrightarrow
(\eta_1,\eta_2)
=
(\eta_{1i},\eta_{2i}),
\quad
1\le i\le 4,
$$
with
$$
(\eta_{11},\eta_{21})
=
(\xi,\xi),
\quad
(\eta_{12},\eta_{22})
=
(\xi,-\xi),
\quad
(\eta_{13},\eta_{23})
=
(-\xi,\xi),
\quad
(\eta_{14},\eta_{24})
=
(\xi/3,\xi/3).
$$
Moreover, taking the second derivative we obtain
\begin{align*}
\p_{\eta_1}^2\Phi
=
3(\xi-\eta_2)(-2)
=
-6(\xi-\eta_2)
\qquad
\p_{\eta_2}^2\Phi
=
3(\xi-\eta_1)(-2)
=
-6(\xi-\eta_1).
\end{align*}

As for the symbol $\Tb_1$ we have the following estimate. From the proof, we will see that $\Tb_1$ is a symbol of order $0$ and $\partial_1\Tb_1$ is a symbol of order $-2$, where the derivative is taken for the variable of the highest frequency. 
\begin{proposition}
\label{Prop_symb_T1}
For $j_1,j_2,j_3\in\Z$ with $j_1\ge j_2\ge j_3$, we have the following estimate 
\begin{align}
\label{T1-est1}
\|\Tb_{1}(\eta_1,\eta_2, \eta_{3})
\psi_{j_1}(\eta_1)\psi_{j_2}(\eta_2)\psi_{j_{3}}(\eta_{3})\|_{S^{\infty}}
\lesssim 
2^{\max\{2j_1,0\}}.
\end{align}
For $\ell=1,2,3$,
\begin{align}
\label{T1-est2}
\|
\p_{\eta_\ell}\Tb_{1}(\eta_1,\eta_2, \eta_{3})
\psi_{j_1}(\eta_1)\psi_{j_2}(\eta_2)\psi_{j_{3}}(\eta_{3})\|_{S^{\infty}}
\lesssim
2^{j_1}.
\end{align}
\end{proposition}

\section{Energy Inequality}
\label{priori-energy-est}
In this section, we prove Lemma \ref{energyLem}.

\subsection{Sobolev energy estimate}
First, we rewrite the equation \eqref{qmkdv:eqn} as following:

\begin{align}
\label{eqn:rewrite}
\vp_t
+
3\vp^2\px \vp
+
\px(c_1(\vp)\px(c_1(\vp)\px\vp))
=
0,
\end{align}
where
\begin{equation}
\label{def-c1}
[c_1(\vp)]^2
\doteq
[c(\vp)]^2+1.
\end{equation}
Notice that $\|\vp(t)\|_{L^2}$ is conserved. The higher-order energy estimate is from \cite{Mie17} by constructing the weighted $L^2$ norms. We include the main steps here for completeness. 

\subsubsection{Higher order energy estimate}
Define the kth-order weighted derivative as $\vp^{(k)}=(c_1(\vp)\px)^k\vp$, and the weight as $\phi_k=c_1^{-\frac{k-1}{3}}$.
The equation of $\vp^{(k)}$ is
\begin{align}
\label{phik-eqn-2}
\phi_k\p_t(\phi_k\vp^{(k)})
+
3\vp^2\phi_k\px(\phi_k\vp^{(k)})
+
\px(\phi_kc_1\px)^2\vp^{(k)}
=
G_k\px\vp^{(k)}
+
H_k,
\end{align}
where
\begin{equation}
\label{def-Gk}
G_k
\doteq
\phi_k^2
g_k
+
\left
[
(\px\phi_k)c_1[\px(\phi_kc_1)]
+
\phi_k[\px(c_1^2\px\phi_k)]
+
\phi_k(\px\phi_k)
c_1(\px c_1)
\right],
\end{equation}
\begin{equation}
\label{def-Hk}
H_k
\doteq
\phi_k^2
h_k
+
(\p_t\phi_k)
\phi_k
\vp^{(k)}
+
3\vp^2(\px\phi_k)\phi_k\vp^{(k)},
\end{equation}
with
\begin{align*}
g_{k}
=&
(k-1)
c_1'\vp_1^2
+
\frac{(k+1)(k-2)}{2}
\px
[c_1'\vp_1]
-
(k+1)(k-2)
\frac{c_{1}{'}^2}{c_1^2}\vp_1^2
+
\left(
\frac{c_1'^2}{c_1^2}
\vp_1^2
-
\frac{c_1'}{c_1}
\vp_2
\right),
\quad
k
\ge
3,
\end{align*}
$$
h_1
=
0,
\quad
\quad
h_2
=
\left(
\frac{6\vp c_1'}{c_1^2}
-
\frac{6}{c_1}
\right)
\vp_1^3
-
\frac{9\vp^2}{c_1}
\vp_1\vp_2
=
\left(
\frac{6\vp c_1'}{c_1^2}
-
\frac{6}{c_1}
\right)
\vp_1^3
-
\frac{6}{c_1}
\vp
\vp_1\vp_2,
$$
\begin{align}
\label{hk+1}
h_{k+1}
\doteq&
c_1\px h_k
-
\frac{6\vp}{c_1}
\vp_1\vp_{k+1}
+
\left[
\frac{c_1'^2}{c_1^3}
f_k
\vp_1^2
-
\frac{c_1'}{c_1^2}
g_k\vp_1
-
\frac{c_1'}{c_1^2}
\px\vp_2
+
\px
\left(
g_k-\frac{c_1'}{c_1}f_k\vp_1
\right)
\right]
\vp_{k+1},
\quad
k\ge 2,
\nonumber
\end{align}
and
$$
f_{k}
=
(k-1)
c_1'\vp_1,
\quad
k
\ge
2.
$$
Multiplying \eqref{phik-eqn-2} by $2\vp^{(k)}$ and 
integrating by parts, we get
\begin{align*}
\frac{d}{dt}
\int_\R
(\phi_k\vp^{(k)})^2
dx
+
\int_\R
3\vp^2
\px(\phi_k\vp^{(k)})^2
dx
+
2\int_\R
[\px(\phi_kc_1\px)^2\vp^{(k)}
]
\vp^{(k)}
dx
=
2\int_\R
G_k
(\px\vp^{(k)})
\vp^{(k)}
dx
+
2\int_\R
H_k
\vp^{(k)}
dx.
\end{align*}
Also, we have 
\begin{align*}
\int_\R
3\vp^2
\px(\phi_k\vp^{(k)})^2
dx
=
-
\int_\R
3\px (\vp^2)
(\phi_k\vp^{(k)})^2
dx,
\end{align*}
\begin{align*}
2\int_\R
[
\px(\phi_kc_1\px)^2\vp^{(k)}
]
\vp^{(k)}
dx
=
-
2\int_\R
\phi_kc_1\px(\phi_kc_1(\px\vp^{(k)}))
\cdot 
(\px\vp^{(k)})
dx
=
-\int
\px
[\phi_kc_1(\px\vp^{(k)})]^2
dx
=
0.
\end{align*}
\begin{align*}
2\int_\R
G_k
(\px\vp^{(k)})
\vp^{(k)}
dx
=
\int_\R
G_k
\cdot
\px(\vp^{(k)})^2 
dx
=
-
\int_\R
\px G_k
\cdot (\vp^{(k)})^2
dx.
\end{align*}
Thus, we get
\begin{align}
\frac{d}{dt}
\|\phi_k\vp^{(k)}\|_{L^2_x}^2
\le
{
3
}
\|\px(\vp^2)\|_{L^\infty}
\|\phi_k\vp^{(k)}\|_{L^2_x}^2
+
\|\px G_k\|_{L^\infty}
\|\vp^{(k)}\|_{L^2_x}^2
+
2
\|H_k\|_{L^2_x}
\|\vp^{(k)}\|_{L^2_x}.
\end{align}
From the definitions of $g_k$ and $G_k$, we obtain that for each $(t, x)$,
\[
|g_k|\lesssim |\vp|_1|\px\vp|_1,
\]
\[
|\px g_k|\lesssim |\px\vp|^2+|\vp\px^2\vp|+|\vp(\px\vp)^3|+|\vp^2\px\vp\px^2\vp|+|\px\vp\px^2\vp|+|\vp\px^3\vp|,
\]
\[
|\px G_k|\lesssim |\px g_k|+|\px\vp g_k|+|\px\vp|^2+|\px\vp\px^2\vp|\lesssim |\px \vp|_2|\vp|_2 F(|\vp|_2),
\]
where $F$ is a positive polynomial. 
Also, we have
\[
\|H_k\|_{L^2}\lesssim \|\vp\|_{H^k}\|F(|\vp|_2) |\px\vp|_2|\vp|_2 \|_{L^\infty}.
\]
Notice that under the smallness assumption, the weighted norm is equivalent to the Sobolev norm, that is for $k\geq 1$,
\[
\|\vp\|_{L^2}+\|\phi_k\vp^{(k)}\|_{L^2}\approx \|\vp\|_{H^k}.
\]
Therefore, the energy inequality can be written as 
\begin{align}\label{HFenergy}
\|\vp(t)\|_{H^k}^2\lesssim \|\vp(0)\|_{H^k}^2+ \int_0^t F(\|\vp\|_{W^{2,\infty}})\||\vp|_2|\px\vp|_2\|_{L^\infty}\|\vp(\tau)\|_{H^k}^2d\tau,
\end{align}
where $F$ is a positive polynomial.

\subsubsection{Low frequency Sobolev estimate}
We will estimate $\|\px^{-1}\vp\|_{L^2}$. First, $\px^{-1}\vp$ satisfies the equation
\begin{align*}
\px^{-1}\vp_t+\px^3(\px^{-1}\vp)
+
\vp^2\px(\px^{-1}\vp)+c(\vp)\px[c(\vp)\px^2(\px^{-1}\vp)]=0.
\end{align*}
Multiplying $\px^{-1}\vp$ and taking integral,
\begin{align}
\nonumber
\frac{d}{dt}\int_{\R} \frac12|\px^{-1}\vp|^2dx=~&\int_{\R}\vp\vp_x (\px^{-1}\vp)^2+c(\vp)\px^2(\px^{-1}\vp)\px[c(\vp)\px^{-1}\vp]dx
\\
\nonumber
=~&\int_{\R}\vp\vp_x (\px^{-1}\vp)^2+\{\px^2[c(\vp)\px^{-1}\vp]-\px^2c(\vp)\px^{-1}\vp-2\px c(\vp)\vp\}\px[c(\vp)\px^{-1}\vp]dx\\\nonumber
=~&\int_{\R}\vp\vp_x (\px^{-1}\vp)^2-
[\px(\px\vp c'(\vp))\px^{-1}\vp+2\px\vp c'(\vp)\vp]
[c'(\vp)\px\vp\px^{-1}\vp+c(\vp)\vp] dx
\\
\label{LFenergy}
\lesssim~& \||\vp|_2|\vp_x|_2\|_{L^\infty}(\|\px^{-1}\vp\|_{L^2}^2+\|\vp\|_{L^2}^2).
\end{align}
%

\subsection{Scaling vector field estimate}
\label{sec:Scaling}
We first estimate $\|S\vp\|_{L^2}$, which can be achieved by acting the operator $S$ on the equation and then do energy estimate.  Then we estimate $\|\px^{-1}S\vp\|_{L^2}$. In the last part of this section, we estimate $\|\xi\partial_\xi\hat h\|_{L^2}$ and $\|\partial_\xi\hat h\|_{L^2}$, which are approximately equivalent to $\|S\vp\|_{L^2}$ and $\|\px^{-1}S\vp\|_{L^2}$ if ignoring a small nonlinear error. 

\subsubsection{$\|S\vp\|_{L^2}$ estimate}

To estimate $\|S\vp\|_{L^2}$, first we summarize some communicator identities for the operator $\S$ defined in \eqref{scaling-operator}.
By direct calculation, we have the following commutator equalities.
\begin{equation}
\label{commu-t-1}
[\S,\p_t]\vp
=
3t\p_t^2\vp
-
\p_t(3t\p_t\vp)
=
-3\p_t\vp,
\end{equation}
\begin{equation}
\label{commu-x-1}
[\S,\p_x]\vp
=
x\p_x^2\vp
-
\p_x(x\p_x\vp)
=
-\p_x\vp,
\end{equation}
\begin{equation}
\label{commu-x-3}
[\S,\p_x^3]\vp
=
x\p_x^4\vp
-
\p_x^3(x\p_x\vp)
=
-3\p_x^3\vp,
\end{equation}
\begin{align}
\label{commu-c-x}
[\S,c(\vp)\px]u
=
\left[
\S\vp
c'(\vp)
-
c(\vp)
\right]\px u,
\end{align}
\begin{equation}
\label{commutators-2}
3\px(\vp^2 S\vp)-\S\px(\vp^3)
=
6\vp\vp_x\S\vp+3\vp^2\partial_x\S\vp-3\vp^2\S\partial_x\vp-6\vp \S\vp\vp_x=3\vp^2 \vp_x,
\end{equation}
and
\begin{align*}
&\px(c(\vp)\px(c(\vp)\px [\S\vp]))-\S\px(c(\vp)\px(c(\vp)\px \vp))\\
=~&\px\bigg(c(\vp)\px([c(\vp)\px,\S] \vp)\bigg)+\px\bigg([c(\vp)\px,\S]c(\vp)\px \vp\bigg)+[\px,\S]c(\vp)\px(c(\vp)\px \vp)\\
=~&\px\Bigg(c(\vp)\px\bigg(\big[c(\vp)-\S\vp c'(\vp)\big]\px\vp\bigg)\Bigg)
+\px\Bigg(\big[c(\vp)-\S\vp c'(\vp)\big]\px\big(c(\vp)\px \vp\big)\Bigg)
+\px(c(\vp)\px(c(\vp)\px \vp))\\
=~&-\px\Bigg(c(\vp)\px\bigg(\big[c'(\vp)\S\vp\big]\px\vp\bigg)\Bigg)
-\px\Bigg(\S\vp c'(\vp)\px\left(c(\vp)\px\vp\right)\Bigg)
+3\px(c(\vp)\px(c(\vp)\px \vp)).
\end{align*}
Therefore $\S\vp$ satisfies following equation:
\begin{align}\nonumber
&\p_t[\S\vp]+\p_x^3[\S\vp]+3\px(\vp^2\S\vp)+\px(c(\vp)\px(c(\vp)\px [\S\vp]))\\\nonumber
=~&[\p_t, \S]\vp+[\p_x^3, \S]\vp+3\px(\vp^2 \S\vp)-\S\px(\vp^3)+\px(c(\vp)\px(c(\vp)\px [\S\vp]))-\S\px(c(\vp)\px(c(\vp)\px \vp))\\\nonumber
=~&3\partial_t\vp+3\px^3\vp+\px(\vp^3)+3\px(c(\vp)\px(c(\vp)\px \vp))-\px\Bigg(c(\vp)\px\bigg(\big[c'(\vp)\S\vp\big]\px\vp\bigg)\Bigg)\\\nonumber
&-\px\Bigg(\S\vp c'(\vp)\px\left(c(\vp)\px\vp\right)\Bigg)\\\nonumber
=~&-2\px(\vp^3)-\px\Bigg(c(\vp)\px\bigg(\big[c'(\vp)\S\vp\big]\px\vp\bigg)\Bigg)-\px\Bigg(\S\vp c'(\vp)\px\left(c(\vp)\px\vp\right)\Bigg)\\\label{eqnS}
=~&b_2\px^2[\S\vp]+b_1\px[\S\vp]+b_0[\S\vp]+F(\vp),
\end{align}
where
\begin{equation*}
b_2
=
b_2(\vp)
\doteq
-\px[c(\vp)]\cdot c(\vp)
=
-
\frac12
\px[c^2(\vp)],
\end{equation*}
\begin{equation*}
b_1
=
b_1(\vp)
\doteq
-cc''\vp_x^2-\px(cc'\vp_x)-c'\px(c\vp_x),
\end{equation*}
\begin{equation*}
b_0
=
b_0(\vp)
\doteq
-\px\left(\px^2[c(\vp)]\cdot c(\vp)+c'(\vp)\px\big(c(\vp)\px \vp\big)\right),
\end{equation*}
and
\begin{equation*}
F(\vp)
\doteq
-2\px(\vp^3).
\end{equation*}

\noindent{\bf  $\|\S\vp\|_{L^2}$ estimate.} We want to find a weight function $w(\vp)$, such that the equation for $w(\vp)\S\vp$ does not have second order term on the right. Multiply $w(\vp)$ to \eqref{eqnS},
\begin{align}
\nonumber
&\p_t[w\S\vp]+\px^3[w\S\vp]+{3\vp^2\px(w\S\vp)}+\px\Big(c(\vp)\px\big(c(\vp)\px[w\S\vp]\big)\Big)\\\nonumber
=~&
w'\vp_t\S\vp+3w'\vp_x[\S\vp]_{xx}+3\px(w'\vp_x)[\S\vp]_{x}+\px^2(w'\vp_x)[\S\vp]+3\vp^2w'\vp_x\S\vp
{-3\px(\varphi^2)w\S\varphi}
\\
\nonumber
&
+\px\Big(c(\vp)\px\big(c(\vp)w'\px\vp\S\vp\big)\Big)
+\px\Big(c(\vp)w'\px\vp\big(c(\vp)\px[\S\vp]\big)\Big)+w'\px\vp \Big(c(\vp)\px\big(c(\vp)\px[\S\vp]\big)\Big)
\\
\nonumber
&+wb_2\px^2[\S\vp]+wb_1\px[\S\vp]+wb_0[\S\vp]+wF(\vp)\\\nonumber
=~&
wF(\vp)
+
w'\S\vp\left[\vp_t+\px^3\vp+3\vp^2\vp_x+\px\Big(c(\vp)\px\big(c(\vp)\px\vp\big)\Big)\right]
+w(b_0
{-3\px(\varphi^2)}
)
\S\vp
\\
\nonumber
&+\px(\S\vp)\left[3\px(w'\vp_x)+2\px\Big(c^2w'\px\vp\Big)+\Big(c(\vp)\px\big(c(\vp)w'\px\vp\big)\Big)+w'\px\vp \Big(c(\vp)\px\big(c(\vp)\big)\Big)+wb_1\right]
\\
\nonumber
&+\px^2(\S\vp)[3w'\px\vp+3c^2w'\px\vp+wb_2]
\\
\label{Seq}
=~&wF(\vp)+a_0^0(\vp)\S\vp+a_1^0(\vp)\px\S\vp+a_2^0(\vp)\px^2\S\vp,
\end{align}
where
\begin{align*}
a_0^0(\vp)
&=
wb_0(\vp)
+
3w
\px(\vp^2)
=
-w\px
\left(
\px^2[c(\vp)]\cdot c(\vp)+c'(\vp)\px\big(c(\vp)\px \vp\big)
{+3\vp^2}
\right),
\\
a_1^0(\vp)&=3\px(w'\vp_x)+2\px\Big(c^2w'\px\vp\Big)+\Big(c(\vp)\px\big(c(\vp)w'\px\vp\big)\Big)+w'\px\vp \Big(c(\vp)\px\big(c(\vp)\big)\Big)+wb_1,
\\
a_2^0(\vp)&=3w'\px\vp+3c^2w'\px\vp+wb_2.
\end{align*}
We will take $w$ such that 
\[
a_2^0(\vp)=3w'\px\vp+3c^2w'\px\vp+wb_2=0,
\]
which leads to 
{
\[
3w'c_1^2\px\vp=\frac12w\px(c_1^2)\Longrightarrow \frac{w'}{w}=\frac16 \frac{2c_1c_1'}{c_1^2}=\frac{c_1'}{3c_1} \Longrightarrow w=c_1^{1/3}.
\]}
Notice that 
\[
|a_0^0(\vp)(t,x)|\lesssim |\vp(t,x)|_{2}|\vp_x(t,x)|_2,\quad |a_1^0(\vp)(t,x)|\lesssim |\vp(t,x)|_{{1}}|\vp_x(t,x)|_1.
\]
Multiplying \eqref{Seq} by $w\S\vp$ and taking integral, we obtain
\begin{align}\nonumber
\frac{d}{dt}\int_{\R} |w\S\vp(t,x)|^2 dx&\lesssim \||\vp|_2|\px\vp|_2\|_{L^\infty} (\|w\S\vp(t)\|_{L^2}^2+\|\vp(t)\|_{L^2}^2)\\\label{S0_est}
&\lesssim \||\vp|_2|\px\vp|_2\|_{L^\infty} E[\vp](t).
\end{align}

\subsubsection{$\|\px^{-1}\S\vp\|_{L^2}$ estimate.}
From \eqref{eqnS}, we have
\begin{multline}\label{xh}
\p_t[\px^{-1}\S\vp]+\px^3[\px^{-1}\S\vp]+3\vp^2\S\vp+\px\left[c(\vp)\px(c(\vp)\px[\px^{-1}\S\vp])\right]-c'(\vp)\px\vp\px(c(\vp)\S\vp)\\
-c(\vp)\px(c'(\vp)\px\vp\S\vp)=-4\vp^3-\Bigg(c(\vp)\px\bigg(\big[c'(\vp)\S\vp\big]\px\vp\bigg)\Bigg)-\Bigg(\S\vp c'(\vp)\px\left(c(\vp)\px\vp\right)\Bigg).
\end{multline}
Multiply this equation by $\px^{-1}\S\vp$ and take integral,
\begin{align*}
\frac{d}{dt}\int_{\R}|\px^{-1}\S\vp|^2dx\lesssim \||\vp|_1|\px\vp|_1\|_{L^\infty}\|\px^{-1}\S\vp\|_{H^1}^2+\|\vp^3\|_{L^2}\|\px^{-1}\S\vp\|_{L^2}.
\end{align*}
Notice that
\begin{align*}
\|\vp^3\|_{L^2}^2=\int_{\R} |\vp(t, x)|^6 dx=\left|\int_{\R}\px[\vp^5]\px^{-1}\vp dx\right|\lesssim \|\px^{-1}\vp\|_{L^2}\|\vp\vp_x\|_{L^\infty}\|\vp^3\|_{L^2},
\end{align*}
which leads to 
\begin{equation}\label{phi3est}
\|\vp^3\|_{L^2}\lesssim \|\px^{-1}\vp\|_{L^2}\|\vp\vp_x\|_{L^\infty}.
\end{equation}

Then by Lemma \ref{lem:decay} and the bootstrap assumption, 
\begin{align}\label{S-1_est}
\frac{d}{dt}\int_{\R}|\px^{-1}\S\vp|^2dx&\lesssim \||\vp|_1|\px\vp|_1\|_{L^\infty}\left(\|\px^{-1}\S\vp\|_{H^1}^2+\|\px^{-1}\vp\|_{L^2}\|\px^{-1}\S\vp\|_{L^2}\right).
\end{align}

\subsection{Estimate of $\|\xi\partial_\xi\hat h(t)\|_{L^2}$ and $\|\partial_\xi\hat h(t)\|_{L^2}$}

Using the definition of $S$ and  \eqref{vp-FT-eqn}, we get
\begin{equation}
\label{x-px-h-fre}
\mathcal{F}_x[x\px h](t, \xi)
=-
\p_\xi(\xi \hat{h}(t, \xi))
=
-\big[
\hat{h}(t, \xi)
+
\xi\p_\xi\hat{h}(t, \xi)
\big],
\end{equation}
and by using \eqref{vp-FT-eqn},
\begin{align}
\label{xi-p-xi-h}
\xi
\p_\xi
\hat{h}(t,\xi)
=&
\xi e^{-it\xi^3}
\left[
-3it\xi^2
\hat{\vp}(t,\xi)
+
\p_\xi\hat{\vp}(t,\xi)
\right]
=
e^{-it\xi^3}
\left[
-3it\xi^3
\hat{\vp}(t,\xi)
+
\xi\p_\xi\hat{\vp}(t,\xi)
\right]
\nonumber
\\
=&
e^{-it\xi^3}
\left[
-3t\p_t\hat{\vp}(t,\xi)
-
3t\widehat{\Nc(\vp)}(t,\xi)
-
\widehat{x\px \vp}(t,\xi)
-
\hat{\vp}(t,\xi)
\right]
\nonumber
\\
=&
-
e^{-it\xi^3}
\left[
\widehat{\S\vp}(t,\xi)
+
3t\widehat{\Nc(\vp)}(t,\xi)
+
\hat{\vp}(t,\xi)
\right].
\end{align}
Combining the above estimates, we obtain
\begin{align}\label{h_est1}
\|x\px h(t)\|_{L^2}
\lesssim
\|S\vp\|_{L^2}
+
3t
\|\Nc(\vp)\|_{L^2}
+
\|\vp\|_{L^2}.
\end{align}
So, we need to estimate the $L^2$-norm of $S\vp$.

As for $\|xh\|_{L^2}$, we have
\[
\|xh\|_{L^2}=\|\frac1\xi\xi\partial_\xi\hat h\|_{L^2}=\|\frac1\xi\xi\partial_\xi\hat h\|_{L^2}=\left\|\frac1\xi
\left[
\widehat{S\vp}(t,\xi)
+
3t\widehat{\Nc(\vp)}(t,\xi)
+
\hat{\vp}(t,\xi)
\right]\right\|_{L^2}.
\]
Notice that 
\[
\px^{-1}\Nc(\vp)=\vp^3+c(\vp)\px[c(\vp)\vp_x].
\]
Then
\begin{align}\label{h_est2}
\|xh\|_{L^2}\lesssim \|\px^{-1}\S\vp\|_{L^2}+t\|\vp^3+c(\vp)\px[c(\vp)\vp_x]\|_{L^2}+\|\px^{-1}\vp\|_{L^2}.
\end{align}
By \eqref{phi3est}
\begin{align}
\|xh\|_{L^2}\lesssim \|\px^{-1}\S\vp\|_{L^2}+t\|\px^{-1}\vp\|_{L^2}\|\vp\vp_x\|_{L^\infty}+t\|c(\vp)\px[c(\vp)\vp_x]\|_{L^2}+\|\px^{-1}\vp\|_{L^2}
\end{align}
\begin{align}
\|c(\vp)\px[c(\vp)\vp_x]\|_{L^2}=\|c(\vp)c'(\vp)\vp_x^2+[c(\vp)]^2\vp_{xx}\|_{L^2}\lesssim \|c\|_{C^1}^2\|\vp\vp_x\|_{L^\infty}\|\vp_x\|_{L^2}+\|[c(\vp)]^2\vp_{xx}\|_{L^2}.
\end{align}
\begin{align}
\|[c(\vp)]^2\vp_{xx}\|_{L^2}^2&=\int c^4\vp_{xx}^2 dx=-\int (c^4\vp_{xx})_x\vp_x dx=-\int 4c^3c'\vp_x^2\vp_{xx}+c^4\vp_{xxx}\vp_xdx\\
&\lesssim \|c\|_{C^1}^4\|\vp\vp_x\|_{L^\infty}\|\vp\|_{H^2}^2+\|c\|_{C^0}^2\|\vp\|_{L^2}^2\|\vp_{xxx}\vp_x\|_{L^\infty}.
\end{align}

By combining \eqref{HFenergy}, \eqref{LFenergy}, \eqref{S0_est}, \eqref{S-1_est}, \eqref{h_est1} and \eqref{h_est2},  we obtain \eqref{EnergyIneq}.

\section{Dispersive Estimate}\label{sec:dis}
The following linear dispersive estimate is from 
 Lemma 2.1 in \cite{GPR16}. See also \cite{CGI19} or \cite{HSZ20p}.
\begin{lemma}
\label{dis-lem}
For any $t\ge 1$, $x\in\R$, $\beta\in[0,1]$, $n\in\mathbb{R}$,
we have the linear dispersive estimates
\begin{align}
\label{LocDis-2}
||\p_x|^\beta e^{-t\partial_x^3} \p_x^nh(t, x) |
&\lesssim 
t^{-1/3-\beta/3}
(1+|x/t^{1/3}|)^{-1/4+\beta/2}
(\|\widehat{|\p_x|^nh}(t)\|_{L^\infty}+t^{-1/6}
\|x|\p_x|^nh(t)\|_{L^2}).
\end{align}
\end{lemma}

Now we use this result to prove Lemma \ref{lem:decay}.
\vskip.1in
\noindent
\begin{proof}[Proof of Lemma \ref{lem:decay}] 
By Lemma \ref{dis-lem}, 
\begin{align}\label{eq5.2}
&|\vp(t,x)|\lesssim t^{-1/3}(1+|x/t^{1/3}|)^{-\frac14}(\|\hat h\|_{L^\infty_\xi}+t^{-\frac16}\|\partial_\xi\hat h\|_{L^2_\xi}),\\\nonumber
&|\vp_x(t,x)|\lesssim t^{-2/3}(1+|x/t^{1/3}|)^{\frac14}(\|\hat h\|_{L^\infty_\xi}+t^{-\frac16}\|\partial_\xi\hat h\|_{L^2_\xi}).
\end{align}
Therefore, using bootstrap assumption, we have
\[
|\vp(t,x)\vp_x(t,x)|\lesssim t^{-1}\ve_1^2.
\]
Since 
\begin{align*}
&\||\px|^{1/2}\vp(t,x)\|_{L^\infty}\lesssim t^{-1/2}(\|\hat h\|_{L^\infty_\xi}+t^{-\frac16}\|\partial_\xi\hat h\|_{L^2_\xi}),
\end{align*}
\begin{align*}
&\|\vp_{xxx}(t,x)\|_{L^\infty}=\|[P_{\{|\xi|\leq t^{1/9-\frac23p_0}\}}+P_{\{|\xi|\geq t^{1/9-\frac23p_0}\}}]\vp_{xxx}(t,x)\|_{L^\infty}\\
 \lesssim& t^{-1/2}(\||\xi|^{5/2}\hat h\|_{L^\infty_\xi}+t^{-p_0}\|\xi\partial_\xi\hat h\|_{L^2_\xi}+t^{-\frac16}\||\xi|^{3/2}\hat h\|_{L^2_\xi})+\||\xi|^{3-\frac s2}(|\xi|^{\frac s2}\hat\vp(\xi))|\|_{L^1_\xi(|\xi|>t^{1/9-\frac23p_0})}\\
 \lesssim& t^{-1/2}(\||\xi|^{5/2}\hat h\|_{L^\infty_\xi}+t^{-p_0}\|\xi\partial_\xi\hat h\|_{L^2_\xi}+t^{-\frac16}\||\xi|^{3/2}\hat h\|_{L^2_\xi})+t^{-\frac{s-7}9+\frac23(s-7)p_0}\|\vp\|_{H^s},
\end{align*}
by bootstrap assumption and interpolation,
\[
\||\px|^{1/2}\vp(t,x)\|_{L^\infty}+\|\vp_x(t,x)\|_{L^\infty}+\|\vp_{xx}(t,x)\|_{L^\infty}+\|\vp_{xxx}(t,x)\|_{L^\infty}\lesssim t^{-1/2}\ve_1.
\]
Therefore,
\[
\||\vp|_2|\vp_x|_2\|_{L^\infty}\lesssim t^{-1}\ve_1^2,
\]
which complete the proof of \eqref{dec1}.

\end{proof}

\section{Z-norm Estimate}\label{sec:Znorm}

In this section,  we prove the estimate
\begin{align}
\label{nonDispest}
(|\xi|^\gl + |\xi|^\gh) |\hat{\vp}(t,\xi)| \lesssim \ve_0\qquad \text{for all $\xi \in \R$},
\end{align}
under the bootstrap assumption \eqref{bootstrap-assump}.
To do this, we define the small number $p_1$
\begin{equation}
\label{p1-choice}
p_1
\doteq
10^{-3}.
\end{equation}
Furthermore, we observe that $p_1$ satisfies the inequality:
$$
p_1 \geq \frac{2p_0}{s+1-2\gh},
$$
where $p_0=10^{-4}$ defined by \eqref{p0-choice}, $s\ge 12$ and $\gh=\frac52$. 

\subsection{Large and small frequencies}
\label{sec:largesmall}

When $|\xi|<(t + 1)^{- 2p_0}$, where $p_0$ is defined by \eqref{p0-choice}, Lemma \ref{interpolation} and the bootstrap assumptions give
\begin{align*}
|(|\xi|^{\gl} + |\xi|^{\gh})\hat\vp(t,\xi)|^2&\lesssim |\xi|^{2\gl-1+2} \||\xi|^{-1}\hat\vp\|_{L^2_\xi}(\|\partial_\xi\hat h\|_{L^2_\xi}+\||\xi|^{-1}\hat\vp\|_{L^2_\xi})
\\
&\lesssim |\xi| \|\px^{-1}\vp(t)\|_{L^2}(\|\partial_\xi\hat h\|_{L^2}+\|\px^{-1}\vp(t)\|_{L^2})\\
&\lesssim \ve_0^2.
\end{align*}

\vskip.1in
\noindent
When $|\xi| > (t + 1)^{p_1}$, where $p_1$ is defined in \eqref{p1-choice}, Lemma \ref{interpolation} and the bootstrap assumptions give
\begin{align*}
|(|\xi|^{\gl} + |\xi|^{\gh}) \hat\vp(t,\xi)|^2 &\lesssim \frac{(|\xi|^{\gl} + |\xi|^{\gh})^2}{|\xi|^{s + 1}} \|\vp\|_{H^s} (
\|
\xi\partial_\xi\hat h
\|_{L^2} + \|\vp\|_{L^2})\\
&\lesssim |\xi|^{2\gh - s-1} \ve_0^2 (t + 1)^{2p_0} \\
&\lesssim  \ve_0^2 .
\end{align*}
Thus, we only need to consider the frequency range
\begin{align}
(t + 1)^{- 2p_0}\lesssim |\xi|\lesssim (t + 1)^{p_1}.
\label{freq_range}
\end{align}
In the following, we fix $\xi$ in this range and denote by $\cutoffxi(t, \xi)$ a smooth cut-off function such that
\begin{align}
\begin{split}
&\text{$\cutoffxi(t, \xi) = 1$ on $\left\{(t, \xi) \mid (t+1)^{-2p_0}\le|\xi|\le(t+1)^{p_1}\right\}$,}
\\
&\text{$\cutoffxi(t, \xi)$ is supported on $\left\{(t, \xi) \mid \frac12(t+1)^{-2p_0}\le|\xi|\le(t+1)^{p_1}+1\right\}$,}
\end{split}
\label{cutoff}
\end{align}
and $|\partial_t\cutoffxi|\lesssim (t+1)^{p_0}$.

%
%

\subsection{Modified scattering}
We will use the method of modified scattering and we start with defining a phase correction
\begin{align}
\label{defTheta}
\Theta(t, \xi)
=&
- \frac{\pi \xi }{3\cdot6\xi} [\Tb_1(\xi, \xi, -\xi) + \Tb_1(\xi,-\xi,\xi) + \Tb_1(-\xi, \xi, \xi)] \int_{1}^t \frac{|\hat \vp(\tau, \xi)|^2}{\tau} \diff{\tau}
\\
=&
-\frac{(2\alpha_2\xi^2-1)\pi}{6}
\int_1^t
\frac{
|\hat{h}(\tau, \xi)|^2
}
{\tau}
d\tau.
\nonumber
\end{align}
Letting 
$$ 
\widehat{v}(t, \xi)
=
e^{i\Theta(t, \xi)}
\hat{h}(t, \xi),
$$
and using \eqref{h-FT-fin} we find that 
\begin{align}
\widehat{v}_t(t, \xi)
=
e^{i\Theta(t, \xi)}
[\hat{h}_t(t, \xi)+i\Theta_t(t, \xi)\hat{h}(t, \xi)]
=U(t, \xi)-e^{i\Theta(t, \xi)}e^{-i\xi^3 t}\widehat{\Nc_{\ge 4}(\vp)}(t, \xi),
\end{align}
where
\begin{align}
\label{defU}
U(t, \xi)
=
e^{i\Theta(t, \xi)}
\Big[
&i\xi
\iint
e^{-it\Phi}
{\bf T_1}(\eta_1,\eta_2,\xi-\eta_1-\eta_2)
\hat{h}(t, \eta_1)
\hat{h}(t, \eta_2)
\hat{h}(t, \xi-\eta_1-\eta_2)
d\eta_1
d\eta_2
\\
-&
\frac{i(2\alpha_2\xi^2-1)\pi}{6}
\frac{
|\hat{h}(t, \xi)|^2
\hat{h}(t, \xi)
}
{t}
\Big].
\nonumber
\end{align}
%
In the following, we will prove for each fixed $\xi$ on the support of $\cutoffxi$, the integral
\[
\left| (|\xi|^{\gl}+|\xi|^{\gh})\int_0^\infty U(t, \xi)-e^{i\Theta(t, \xi)}e^{-i\xi^3 t}\widehat{\Nc_{\ge 4}(\vp)}(t, \xi)dt \right| 
\]
is bounded. Assume $|\xi|\in(2^{j-1}, 2^j]$, we will prove the bound of the above integral is uniform in $j$ for different situations. 

%
We first do the dyadic decomposition to the integral. Without loss of generality we assume that $j_1\ge j_2\ge j_3$.
Carrying out a dyadic decomposition, with $h_j=P_j h$, where $P_j$ is the Fourier multiplier with symbol $\psi_j$ defined in \eqref{defpsik}, we have
\begin{align*}
&i\xi
\iint
e^{-it\Phi}
{\bf T_1}(\eta_1,\eta_2,\xi-\eta_1-\eta_2)
\hat{h}(t, \eta_1)
\hat{h}(t, \eta_2)
\hat{h}(t, \xi-\eta_1-\eta_2)
d\eta_1
d\eta_2
\\
=&
\sum\limits_{\P}
i\xi
\iint
e^{-it\Phi}\,{\bf T_1}(\eta_1,\eta_2,\xi-\eta_1-\eta_2)
\hat{h}_{j_1}(t, \eta_1)
\hat{h}_{j_2}(t, \eta_2)
\hat{h}_{j_3}(t, \xi-\eta_1-\eta_2)
d\eta_1
d\eta_2,
\end{align*}
where   $\P$ denotes the set with the indices $j_1,j_2,j_3\in\Z$, such that $j_1\geq j_2\geq j_3$ with possible repetition.
To estimate this integral, we split $\P$ into several cases. 
\begin{align*}
\P_1&=\{(j_1, j_2, j_3)\in\P\mid  j_1-j_3>2\},\\
\P_2&=\{(j_1, j_2, j_3)\in\P\mid  j_1-2\leq j_3\leq j_2\leq j_1\}.
\end{align*}

The following of this section is arranged as follows. The  term $U$ is the cubic nonlinearity, whose estimate is in sections \ref{nonre} to \ref{sptmr}. 
Section \ref{nonre} is the estimate for non-resonant frequencies. When the frequencies are close to resonances, we need to divide the domain further. The domain with a distance away from resonances is discussed in section \ref{clos}. The estimate around the space resonance is in section \ref{spres}, and the space-time resonances in section \ref{sptmr}. The estimate of higher degree terms is in section \ref{higherdegree}.

\subsection{Non-resonant frequencies}\label{nonre}
When $(j_1, j_2, j_3)\in \P_1$, that is $j_1-j_3>2$, we  will estimate
\begin{align}\label{7.12}
\left\|\cutoffxi(|\xi|^\gl+|\xi|^\gh)\xi\int_0^t \iint_{\R^2}   \Tb_1(\eta_1, \eta_2, \xi - \eta_1 - \eta_2) e^{i  \tau \Phi} \hat h_{j_1}(\eta_1)\hat h_{j_2}(\eta_2) \hat h_{j_3}(\xi-\eta_1-\eta_2) \diff{\eta_1} \diff{\eta_2}\diff\tau \right\|_{L^\infty_\xi}.
\end{align}

We rewrite the integral as
\begin{align*}
&  \iint_{\R^2}    \Tb_1(\eta_1, \eta_2, \xi - \eta_1 - \eta_2) e^{it\Phi(\xi,\eta_1,\eta_2)} \hat h_{j_1}(\eta_1)\hat h_{j_2}(\eta_2) \hat h_{j_3}(\xi-\eta_1-\eta_2) \diff{\eta_1} \diff{\eta_2}\\
 =~&   \iint_{\R^2}    \frac{\Tb_1(\eta_1, \eta_2, \xi - \eta_1 - \eta_2)}{it\partial_{\eta_1}\Phi(\xi,\eta_1,\eta_2)}\partial_{\eta_1} e^{it\Phi(\xi,\eta_1,\eta_2)} \hat h_{j_1}(\eta_1)\hat h_{j_2}(\eta_2) \hat h_{j_3}(\xi-\eta_1-\eta_2) \diff{\eta_1} \diff{\eta_2}\\
 =~& -{\rm I}_1-{\rm I}_2-{\rm I}_3,
\end{align*}
where
\begin{align}
\label{non-reson-W1}
{\rm I}_1(t, \xi) & =  \iint_{\R^2} \partial_{\eta_1}\left[ \frac{\Tb_1(\eta_1, \eta_2, \xi - \eta_1 - \eta_2)}{it\partial_{\eta_1}\Phi(\xi,\eta_1,\eta_2)}\right] e^{it\Phi(\xi,\eta_1,\eta_2)}  \hat h_{j_1}(\eta_1)\hat h_{j_2}(\eta_2) \hat h_{j_3}(\xi-\eta_1-\eta_2) \diff{\eta_1} \diff{\eta_2},
\\
\label{non-reson-W2}
{\rm I}_2(t, \xi) & =  \iint_{\R^2} \left[ \frac{\Tb_1(\eta_1, \eta_2, \xi - \eta_1 - \eta_2)}{it\partial_{\eta_1}\Phi(\xi,\eta_1,\eta_2)}\right] e^{it\Phi(\xi,\eta_1,\eta_2)} \hat h_{j_1}(\eta_1)\hat h_{j_2}(\eta_2) \partial_{\eta_1}\hat h_{j_3}(\xi-\eta_1-\eta_2) \diff{\eta_1} \diff{\eta_2},
\\
\label{non-reson-W3}
{\rm I}_3(t, \xi) & =  \iint_{\R^2} \left[ \frac{\Tb_1(\eta_1, \eta_2, \xi - \eta_1 - \eta_2)}{it\partial_{\eta_1}\Phi(\xi,\eta_1,\eta_2)}\right] e^{it\Phi(\xi,\eta_1,\eta_2)} \partial_{\eta_1} \hat h_{j_1}(\eta_1)\hat h_{j_2}(\eta_2) \hat h_{j_3}(\xi-\eta_1-\eta_2) \diff{\eta_1} \diff{\eta_2}.
\end{align}
In the first integral $W_1$, the symbol can be written as two parts $m_1$ and $m_2$ as below
\begin{align*}
\partial_{\eta_1}\left[ \frac{\Tb_1(\eta_1, \eta_2, \xi - \eta_1 - \eta_2)}{\partial_{\eta_1}\Phi(\xi,\eta_1,\eta_2)}\right] 
=&
\frac{\partial_{\eta_1}\Tb_1(\eta_1, \eta_2, \xi - \eta_1 - \eta_2)}{\partial_{\eta_1}\Phi(\xi,\eta_1,\eta_2)}
-
\frac{\Tb_1(\eta_1, \eta_2, \xi - \eta_1 - \eta_2)\partial_{\eta_1}^2\Phi(\xi,\eta_1,\eta_2)}{[\partial_{\eta_1}\Phi(\xi,\eta_1,\eta_2)]^2}
\\
\doteq&m_1(\eta_1, \eta_2, \xi - \eta_1 - \eta_2)+m_2(\eta_1, \eta_2, \xi - \eta_1 - \eta_2),
\end{align*}
where
\begin{align*}
m_1(\eta_1, \eta_2, \xi - \eta_1 - \eta_2)
=
\frac{\partial_{\eta_1}\Tb_1(\eta_1, \eta_2, \xi - \eta_1 - \eta_2)}{\partial_{\eta_1}\Phi(\xi,\eta_1,\eta_2)}
=
\frac{\alpha_2(2\eta_1+\eta_2-\xi)}{9(\xi-\eta_2)(\xi-2\eta_1-\eta_2)}
=
-\frac{\alpha_2}{9(\xi-\eta_2)},
\end{align*}
and 
\begin{align*}
m_2(\eta_1, \eta_2, \xi - \eta_1 - \eta_2)
&=
-
\frac{\Tb_1(\eta_1, \eta_2, \xi - \eta_1 - \eta_2)\partial_{\eta_1}^2\Phi(\xi,\eta_1,\eta_2)}{[\partial_{\eta_1}\Phi(\xi,\eta_1,\eta_2)]^2}\\
&=
\left[\frac29\alpha_2(\eta_1^2+\eta_2^2+\xi^2+\eta_1\eta_2-\eta_1\xi-\eta_2\xi)-\frac23\right][(\xi-\eta_2)(2\eta_1+\eta_2-\xi)^2]^{-1}.
\end{align*}
The corresponding two integrals are denoted as
\begin{align}
{\rm I}_{11}(t, \xi) 
=&
\frac{1}{it}
\iint_{\R^2} m_1(\eta_1, \eta_2, \xi - \eta_1 - \eta_2)e^{it\Phi(\xi,\eta_1,\eta_2)}  \hat h_{j_1}(\eta_1)\hat h_{j_2}(\eta_2) \hat h_{j_3}(\xi-\eta_1-\eta_2) \diff{\eta_1} \diff{\eta_2},
\end{align}
and
\begin{align}
{\rm I}_{12}(t, \xi) 
=&
\frac{1}{it}
\iint_{\R^2} m_2(\eta_1, \eta_2, \xi - \eta_1 - \eta_2)e^{it\Phi(\xi,\eta_1,\eta_2)}  \hat h_{j_1}(\eta_1)\hat h_{j_2}(\eta_2) \hat h_{j_3}(\xi-\eta_1-\eta_2) \diff{\eta_1} \diff{\eta_2}.
\end{align}
We notice that
\begin{align*}
2^{j}
\lesssim
|\xi|
=&
|\eta_1+\eta_2+\eta_3|
\le
|\eta_1|+|\eta_2|+|\eta_3|
\le
3
\cdot
2^{j_1},
\\
|\xi-\eta_2|
=&
|\xi-\eta_1-\eta_2
+
\eta_1
|
=
|\eta_3+\eta_1|
\ge
|\eta_1|
-
|\eta_3|
=
2^{j_1}
(1-2^{j_3-j_1})
\ge
\frac34
\cdot
2^{j_1},
\\
|\xi-2\eta_1-\eta_2|
=&
|\xi-\eta_1-\eta_2
-
\eta_1
|
=
|\eta_3-\eta_1|
\ge
|\eta_1|
-
|\eta_3|
=
2^{j_1}
(1-2^{j_3-j_1})
\ge
\frac34
\cdot
2^{j_1}.
\end{align*}
Also, applying Lemma \ref{multilinear}, we get
\begin{align}\nonumber
\|(|\xi|^\gl+|\xi|^\gh)\xi {\rm I}_{11}(t, \xi)\|_{L^\infty_\xi}
& \lesssim 
(1+t)^{-1}(1+2^{\gh j})
2^j2^{-j_1}
\|\vp_{j_1}\|_{L^\infty}\|\vp_{j_2}\|_{L^2}\|\vp_{j_3}\|_{L^2}
\\\label{est_I11}
& \lesssim 
(1+t)^{-1+\gh p_1+4p_0}
\|2^{2j_1}\vp_{j_1}\|_{L^\infty}\|\vp_{j_2}\|_{L^2}\|\vp_{j_3}\|_{L^2}.
\end{align}
Similarly, for ${\rm I}_{12}$, it satisfies
\begin{align}\nonumber
&\|(|\xi|^\gl+|\xi|^\gh)\xi {\rm I}_{12}(t, \xi)\|_{L^\infty_\xi}
\lesssim 
(1+t)^{-1}
2^j(1+2^{\gh j})2^{-3j_1}\|\vp_{j_1}\|_{L^\infty}\|\vp_{j_2}\|_{L^2}\|\vp_{j_3}\|_{L^2}
\\\label{est_I12}
\lesssim~ &
(1+t)^{-1}(1+2^{\gh j})2^{-4j}
\|2^{2j_1}\vp_{j_1}\|_{L^\infty}\|\vp_{j_2}\|_{L^2}\|2^{-j_3}\vp_{j_3}\|_{L^2}\\
\lesssim~&(1+t)^{-1+\gh p_1+8p_0}\|2^{2j_1}\vp_{j_1}\|_{L^\infty}\|\vp_{j_2}\|_{L^2}\|2^{-j_3}\vp_{j_3}\|_{L^2}.
\end{align}
The symbol for ${\rm I}_2$ and ${\rm I}_3$ is 
$$
\frac{\Tb_1(\eta_1, \eta_2, \xi - \eta_1 - \eta_2)}{\partial_{\eta_1}\Phi(\xi,\eta_1,\eta_2)}
=
\frac{
\frac{\alpha_2}3(\eta_1^2+\eta_2^2+\xi^2+\eta_1\eta_2-\eta_1\xi-\eta_2\xi)-1
}{
3(\xi-\eta_2)(\xi-2\eta_1-\eta_2)
},
$$
which is bounded as follows
$$
\left\|
\frac{\Tb_1(\eta_1, \eta_2, \xi - \eta_1 - \eta_2)}{\partial_{\eta_1}\Phi(\xi,\eta_1,\eta_2)}
\right\|_{S^\infty}
\lesssim
(2^{2j_1}+1)2^{-2j_1}
=
1+2^{-2j_1}.
$$
This estimate implies that 
\begin{align}\nonumber
\|(|\xi|^\gl+|\xi|^\gh)\xi {\rm I}_{2}(t, \xi)\|_{L^\infty_\xi} 
\lesssim ~&
(1+t)^{-1}
2^j(1+2^{\gh j})(1+2^{-2j_1})
\|\vp_{j_1}\|_{L^2}\|\vp_{j_2}\|_{L^\infty}\|\partial_{\xi}
\hat{h}_{j_3}\|_{L^2}\\\label{est_I2}
\lesssim~&(1+t)^{-1}2^j(1+2^{\gh j})(1+2^{-2j})\|\vp_{j_1}\|_{L^\infty}\|\vp_{j_2}\|_{L^2}\|\partial_{\xi}
\hat{h}_{j_3}\|_{L^2}\\
\lesssim~&(1+t)^{-1}(2^{-3j}+2^{\gh j})\|2^{2j_1}\vp_{j_1}\|_{L^\infty}\|\vp_{j_2}\|_{L^2}\|\partial_{\xi}
\hat{h}_{j_3}\|_{L^2}\\
\lesssim~&(1+t)^{-1+6p_0+\gh p_1}\|2^{2j_1}\vp_{j_1}\|_{L^\infty}\|\vp_{j_2}\|_{L^2}\|\partial_{\xi}
\hat{h}_{j_3}\|_{L^2}.
\end{align}
The similar estimate holds for $\rm I_3$:
\begin{align}\nonumber
\|(|\xi|^\gl+|\xi|^\gh)\xi {\rm I}_{3}(t, \xi)\|_{L^\infty_\xi} 
\lesssim~ &
(1+t)^{-1}
2^j(1+2^{\gh j})(1+2^{-2j_1})
\|\vp_{j_3}\|_{L^2}\|\vp_{j_2}\|_{L^\infty}\|\partial_{\xi}
\hat{h}_{j_1}\|_{L^2}
\\\label{est_I3}
\lesssim ~&(1+t)^{-1+(\gh+1) p_1+4p_0}\|2^{-j_3}\vp_{j_3}\|_{L^2}\|2^{j_2}\vp_{j_2}\|_{L^\infty}\|\partial_{\xi}
\hat{h}_{j_1}\|_{L^2}.
\end{align}
Using Lemma \ref{lem:decay}, by taking the summation of \eqref{est_I11}, \eqref{est_I12}, \eqref{est_I2} and \eqref{est_I3}, for $\P_1$ we obtain
\begin{align*}
&
\int_0^\infty\left\|(|\xi|^\gl+|\xi|^\gh)\xi \sum\limits_{
\P_1}\iint_{\R^2}    \Tb_1(\eta_1, \eta_2, \xi - \eta_1 - \eta_2) e^{it\Phi(\eta_1,\eta_2, \xi)} \hat h_{j_1}(\eta_1)\hat h_{j_2}(\eta_2) \hat h_{j_3}(\xi-\eta_1-\eta_2) \diff{\eta_1} \diff{\eta_2}\right\|_{L^\infty_\xi}dt\\
\lesssim~&\int_0^\infty (1+t)^{-1+8p_0+\gh p_1} (\|\px^{-1}\vp(t)\|_{L^2}+\|\vp(t)\|_{L^2}+\|\partial_\xi\hat h(t)\|_{L^2})^2\|\vp(t)\|_{B^{1,2}} dt\\
\lesssim~&\ve_1^3 \int_0^\infty (1+t)^{-1+8p_0+\gh p_1-\frac12} dt\lesssim\ve_1^3\lesssim 
\ve_0.
\end{align*}
where the last line is from the bootstrap assumption \eqref{bootstrap-assump}.

\subsection{Close to Resonances} \label{clos}
When $j_1-j_3\leq 2$, we need further  decomposition to the region close to resonances. To do this, we define 
\[
\upsilon_\pm(\eta)=\left\{
\begin{aligned}
&1,\quad {\text {if } } \pm\eta\geq 0,\\
&0, \quad {\text {if } } \pm\eta<0,
\end{aligned}
\right.
\]
and $\hat h_{j}^\pm(\eta)=\hat h_j(\eta)\upsilon_\pm(\eta)$. 
 For $\iota_1, \iota_2, \iota_3\in \{+,-\}$, we split the integral into several parts
 \begin{align*}
\iint_{\R^2}    \Tb_1(\eta_1, \eta_2, \xi - \eta_1 - \eta_2) e^{it\Phi(\eta_1,\eta_2, \xi)} \hat h_{j_1}^{\iota_1}(\eta_1)\hat h_{j_2}^{\iota_2}(\eta_2) \hat h_{j_3}(\xi-\eta_1-\eta_2)\upsilon_{\iota_3}(\xi) \diff{\eta_1} \diff{\eta_2}.
 \end{align*}
And we will estimate each case in the following.

In the following, we assume that $\xi>0$. The case when $\xi<0$ can be discussed in the similar way.

\noindent{\bf The case of $(-,-)$.} When $\eta_1<0, \eta_2<0, \xi>0$, we have $|\xi-\eta_1-\eta_2|=\xi-\eta_1-\eta_2>0$, thus 
\[
|\partial_{\eta_1}\Phi|\gtrsim 2^{2j_1}.
\]
This is the case away from resonances, whose estimate can be achieved in the same way as in Section~\ref{nonre}.

\noindent{\bf The case of $(+,+)$.} The space-time resonance $(\eta_1, \eta_2)=(\xi,\xi)$ and the space resonance $(\eta_1, \eta_2)=(\xi/3,\xi/3)$ are both in the support of the integrand. By introducing an additional cut-off function $\psi_{j_1-3}(\eta_1+\eta_2-2\xi)$, we write the integral into two parts:
\begin{align*}
\iint_{\R^2} \Tb_1(\eta_1, \eta_2, \xi - \eta_1 - \eta_2) e^{i  t \Phi(\xi,\eta_1,\eta_2)} \hat h_{j_1}^+(\eta_1)\hat h_{j_2}^+(\eta_2) \hat h_{j_3}(\xi-\eta_1-\eta_2)\psi_{j_1-3}(\eta_1+\eta_2-2\xi)\diff\eta_1\diff\eta_2,\\
\iint_{\R^2} \Tb_1(\eta_1, \eta_2, \xi - \eta_1 - \eta_2) e^{i  t \Phi(\xi,\eta_1,\eta_2)} \hat h_{j_1}^+(\eta_1)\hat h_{j_2}^+(\eta_2) \hat h_{j_3}(\xi-\eta_1-\eta_2)[1-\psi_{j_1-3}(\eta_1+\eta_2-2\xi)]\diff\eta_1\diff\eta_2,
\end{align*}
where the first integral contains the space-time resonance $(\eta_1, \eta_2)=(\xi,\xi)$ and the second integral contains the space resonance $(\eta_1, \eta_2)=(\xi/3,\xi/3)$.
Then we estimate these two integrals in the following.  

We consider the partition of unity
\[
\sum\limits_{(k_1,k_2)\in\Z^2}\psi_{k_1}(\eta_1-\xi)\psi_{k_2}(\eta_2-\xi)=1 \quad\text{and}\quad \sum\limits_{(k_3,k_4)\in\Z^2}\psi_{k_3}\bigg(\eta_1-\frac\xi3\bigg)\psi_{k_4}\bigg(\eta_2-\frac\xi3\bigg)=1.
\]
Then we write the integrals with finer dyadic decomposition as
\begin{equation}\label{stRes}
\begin{aligned}
&\iint_{\R^2} \Tb_1(\eta_1, \eta_2, \xi - \eta_1 - \eta_2) e^{i t \Phi(\xi,\eta_1,\eta_2)} \hat h_{j_1}^+(\eta_1)\hat h_{j_2}^+(\eta_2) \hat h_{j_3}(\xi-\eta_1-\eta_2)\psi_{j_1-3}(\eta_1+\eta_2-2\xi)\\
&\qquad\psi_{k_1}(\eta_1-\xi)\psi_{k_2}(\eta_2-\xi)\diff\eta_1\diff\eta_2,
\end{aligned}
\end{equation}
\begin{equation}\label{spRes}
\begin{aligned}
&\iint_{\R^2} \Tb_1(\eta_1, \eta_2, \xi - \eta_1 - \eta_2) e^{i  t \Phi(\xi,\eta_1,\eta_2)} \hat h_{j_1}^+(\eta_1)\hat h_{j_2}^+(\eta_2) \hat h_{j_3}(\xi-\eta_1-\eta_2)[1-\psi_{j_1-3}(\eta_1+\eta_2-2\xi)]\\
&\qquad\psi_{k_3}\bigg(\eta_1-\frac\xi3\bigg)\psi_{k_4}\bigg(\eta_2-\frac\xi3\bigg)\diff\eta_1\diff\eta_2.
\end{aligned}
\end{equation}
In this subsection, we introduce a decay function
\begin{equation}
 \varrho_1(t)=(t+1)^{-0.48},
 \label{def_rho}
 \end{equation}
and restrict our attention to the following two near-resonant cases: 
%
\begin{align*}
&\text{(i) $\max\{k_1, k_2\}\geq \log_2[\varrho_1(t)]$},
\\
&\text{(ii) $\max\{k_3, k_4\}\geq \log_2[\varrho_1(t)]$}.
\end{align*}
The cases of $\max\{k_1,k_2\}< \log_2[\varrho_1(t)]$ or $\max\{k_3,k_4\}< \log_2[\varrho_1(t)]$, which include the resonant frequencies, will be discussed in Section \ref{spres} and Section \ref{sptmr}.

\noindent {\it Near space-time resonance.}
We first estimate \eqref{stRes}.
Since the integrals are symmetric in $\eta_1$ and $\eta_2$, we can assume without loss of generality that $ k_2\geq k_1$.
Using integration by parts, we write the integral as
\begin{equation*}
\begin{aligned}
&
\iint_{\R^2}   \frac{\Tb_1(\eta_1, \eta_2, \xi - \eta_1 - \eta_2)}{3 i  t ( (\xi-\eta_1-\eta_2)^{2}-\eta_1^{2} )}\partial_{\eta_1} e^{i  t \Phi(\xi,\eta_1,\eta_2)} \hat h_{j_1}^+(\eta_1)\hat h_{j_2}^+(\eta_2) \hat h_{j_3}(\xi-\eta_1-\eta_2)
\\
&\qquad \cdot \psi_{j_1-3}(\eta_1+\eta_2-2\xi)\bigg[\psi_{\le k_2}(\eta_1-\xi)\psi_{k_2}(\eta_2-\xi) \bigg] \diff{\eta_1} \diff{\eta_2}
\\
  =~& \frac{V_1+V_2+V_3+V_4}{3i t},
\end{aligned}
\end{equation*}
where
\begin{align*}
V_1(t, \xi)=&\iint_{\R^2}   \partial_{\eta_1}\left[\frac{\Tb_1(\eta_1, \eta_2, \xi - \eta_1 - \eta_2)}{|\eta_1|^{2} - |\xi-\eta_1-\eta_2|^{2}}\right] e^{i  t \Phi(\xi,\eta_1,\eta_2)}  \hat h_{j_1}^+(\eta_1)\hat h_{j_2}^+(\eta_2) \hat h_{j_3}(\xi-\eta_1-\eta_2) \\*
& \hspace{2in} \cdot \psi_{j_1-3}(\eta_1+\eta_2-2\xi)\psi_{\le k_2}(\eta_1-\xi)\psi_{k_2}(\eta_2-\xi) \diff{\eta_1} \diff{\eta_2},
\\[1ex]
   V_2(t, \xi)=& \iint_{\R^2}   \left[\frac{\Tb_1(\eta_1, \eta_2, \xi - \eta_1 - \eta_2)}{|\eta_1|^{2} - |\xi-\eta_1-\eta_2|^{2}}\right] e^{i  t \Phi(\xi,\eta_1,\eta_2)} \partial_{\eta_1} \hat h_{j_1}^+(\eta_1)\hat h_{j_2}^+(\eta_2) \hat h_{j_3}(\xi-\eta_1-\eta_2) 
   \\*
   & \hspace{2in} \cdot \psi_{j_1-3}(\eta_1+\eta_2-2\xi)\psi_{\le k_2}(\eta_1-\xi)\psi_{k_2}(\eta_2-\xi) \diff{\eta_1} \diff{\eta_2},\\[1ex]
  V_3(t, \xi)=& \iint_{\R^2}   \left[\frac{\Tb_1(\eta_1, \eta_2, \xi - \eta_1 - \eta_2)}{|\eta_1|^{2} - |\xi-\eta_1-\eta_2|^{2}}\right] e^{i t \Phi(\xi,\eta_1,\eta_2)}  \hat h_{j_1}^+(\eta_1)\hat h_{j_2}^+(\eta_2)\partial_{\eta_1} \hat h_{j_3}(\xi-\eta_1-\eta_2) 
  \\*
  & \hspace{2in} \cdot \psi_{j_1-3}(\eta_1+\eta_2-2\xi)\psi_{\le k_2}(\eta_1-\xi)\psi_{k_2}(\eta_2-\xi)  \diff{\eta_1} \diff{\eta_2},\\[1ex]
  V_4(t, \xi)=& \iint_{\R^2}   \left[\frac{\Tb_1(\eta_1, \eta_2, \xi - \eta_1 - \eta_2)}{|\eta_1|^{2} - |\xi-\eta_1-\eta_2|^{2}}\right] e^{i  t \Phi(\xi,\eta_1,\eta_2)}  \hat h_{j_1}^+(\eta_1)\hat h_{j_2}^+(\eta_2) \hat h_{j_3}(\xi-\eta_1-\eta_2) 
  \\*
  & \hspace{2in} \cdot \partial_{\eta_1}\bigg[\psi_{j_1-3}(\eta_1+\eta_2-2\xi)\psi_{\le k_2}(\eta_1-\xi)\psi_{k_2}(\eta_2-\xi)\bigg]  \diff{\eta_1} \diff{\eta_2}.
\end{align*}

Then we estimate $V_1$ to $V_4$. We first denote the symbol of $V_1$ by
\begin{align*}
m(\eta_1,\eta_2,\xi-\eta_1-\eta_2)&=\partial_{\eta_1}\left[\frac{\Tb_1(\eta_1, \eta_2, \xi - \eta_1 - \eta_2)}{|\eta_1|^{2} - |\xi-\eta_1-\eta_2|^{2}}\right] \\
&=\frac{\partial_{\eta_1}\Tb_1(\eta_1, \eta_2, \xi - \eta_1 - \eta_2)}{|\eta_1|^{2} - |\xi-\eta_1-\eta_2|^{2}}-2\Tb_1(\eta_1, \eta_2, \xi - \eta_1 - \eta_2)\frac{\eta_1+(\xi-\eta_1-\eta_2)}{(|\eta_1|^{2} - |\xi-\eta_1-\eta_2|^{2})^2},\\
&=\frac{\alpha_2}{3(\xi-\eta_2)}-2\Tb_1(\eta_1, \eta_2, \xi - \eta_1 - \eta_2)\frac{1}{(2\eta_1+\eta_2-\xi)^2(\xi-\eta_2)}.
\end{align*}

%
Then we have the symbol estimate
\begin{align*}
&\Bigg\|m(\eta_1,\eta_2,\eta_3)\psi_{j_1}(\eta_1)\psi_{j_2}(\eta_2)\psi_{j_3}(\eta_3)\upsilon_+(\eta_1)\upsilon_+(\eta_2)\psi_{\le k_2}(-\eta_2-\eta_3)\psi_{k_2}(-\eta_1-\eta_3)\Bigg\|_{S^\infty}\\
\lesssim ~&2^{-k_2}+2^{-2j_1}(2^{2j_1}+1)2^{-k_2}\lesssim 2^{-k_2}(2^{-2j_1}+1).
\end{align*}

Therefore, 
\begin{align*}
\|V_1\|_{L^\infty_\xi}
&\lesssim 2^{-k_2}(2^{-2j_1}+1)\| \vp_{j_1}\|_{L^2} \left\| \vp_{j_2}\right\|_{L^2} \| \vp_{j_3}\|_{L^\infty}.
\end{align*}
So we have
\begin{align}\label{est_V1}
\||\xi|(1+|\xi|^{\gh})\cutoffxi V_1\|_{L^\infty_\xi}
&\lesssim 2^{-k_2}(t+1)^{\gh p_1}\left\|(1+2^{-j_1}) \vp_{j_1}\right\|_{L^2} \left\|(1+2^{-j_2})  \vp_{j_2}\right\|_{L^2} \|\px \vp_{j_3}\|_{L^\infty}.
\end{align}


For the term $V_2$, since its symbol satisfies
\begin{align*}
\left\|\frac{\Tb_1(\eta_1,\eta_2,\eta_3)}{\eta_1^2-\eta_3^2}\psi_{j_1}(\eta_1)\psi_{j_2}(\eta_2)\psi_{j_3}(\eta_3)\upsilon_+(\eta_1)\upsilon_+(\eta_2)\psi_{\leq k_2}(-\eta_2-\eta_3)\psi_{k_2}(-\eta_1-\eta_3)\right\|_{S^\infty}\lesssim (1+2^{2j_1})2^{-j_1-k_2}, 
\end{align*}
considering the bound of $2^{j_1}$, 
\begin{align}\label{est_V2}
&\||\xi|(1+|\xi|^{\gh})\cutoffxi V_2\|_{L^\infty_\xi}\lesssim 2^{-k_2}(t+1)^{(\gh+2) p_1+2p_0}\|\partial_{\eta_1}\hat h_{j_1}\|_{L^2} \left\| \vp_{j_2}\right\|_{L^2} \| \px\vp_{j_3}\|_{L^\infty}.
\end{align}

Similarly we obtain the estimates for $V_3$ and $V_4$ by using the above symbol estimates.
\begin{align}\label{est_V3}
&\||\xi|(1+|\xi|^{\gh})\cutoffxi V_3\|_{L^\infty_\xi}\lesssim 2^{-k_2}(t+1)^{(\gh+2) p_1+2p_0}\|\px \vp_{j_1}\|_{L^\infty} \left\| \vp_{j_2}\right\|_{L^2}\|\partial_{\eta_3}\hat h_{j_3}\|_{L^2},
\end{align}
\begin{align}\nonumber
&\||\xi|(1+|\xi|^{\gh})\cutoffxi V_4\|_{L^\infty_\xi}\\\nonumber
\lesssim ~&2^{-k_2}(t+1)^{(\gh+2) p_1+2p_0}( 2^{-k_2}+2^{-j_1})\|\psi_{\le k_2}(\eta_1-\xi) \hat\vp_{j_1}(\eta_1)\|_{L^2_{\eta_1}L^\infty_\xi} \|\psi_{k_2}(\eta_2-\xi) \hat\vp_{j_2}(\eta_2)\|_{L^2_{\eta_2}L^\infty_\xi} \| \vp_{j_3}\|_{L^\infty}
\\\label{est_V4}
\lesssim~& (2^{-k_2}+2^{-j_1})(t+1)^{2p_0+(\gh+4) p_1} \|\hat\vp_{j_1}\|_{L^\infty_\xi}\|\hat\vp_{j_2}\|_{L^\infty_\xi}\| \px\vp_{j_3}\|_{L^\infty}.
\end{align}

Finally, we combine the estimates from \eqref{est_V1} to \eqref{est_V4} and  take the summation over $\log_2[\varrho_1(t)] \leq k_2 \leq j_1 + 5$. Notice that 
\begin{align}\label{6.25}
\sum_{k_2=\log_2[\varrho_1(t)]}^{j_1+5}2^{-k_2}\lesssim \max\{2^{-j_1-5}, [\varrho_1(t)]^{-1}\}\lesssim (t+1)^{0.48},
\end{align}
and
\begin{align}\label{6.26}
\sum_{k_2=\log_2[\varrho_1(t)]}^{j_1+5}2^{-j_1}\lesssim 2^{-j_1}\log(t+1)\lesssim (t+1)^{2p_0+p_1}.
\end{align}
So we have the estimate
\begin{equation}
\label{eqn829}
\begin{aligned}
& \Bigg\|\cutoffxi\xi (|\xi|^{\gl} + |\xi|^{\gh})\int_0^t  \iint_{\R^2} \Tb_1(\eta_1, \eta_2, \xi - \eta_1 - \eta_2) e^{i t \Phi(\xi,\eta_1,\eta_2)} \hat h_{j_1}^+(\eta_1)\hat h_{j_2}^+(\eta_2) \hat h_{j_3}(\xi-\eta_1-\eta_2) 
\\
& \hspace{1.5in} \cdot\psi_{j_1-3}(\eta_1+\eta_2-2\xi) \bigg[\sum\limits_{\substack{k_2\geq \log_2[\varrho_1(\tau)]\\k_1\le k_2}}\psi_{k_1}(\eta_1-\xi)\psi_{k_2}(\eta_2-\xi) \bigg] \diff{\eta_1} \diff{\eta_2}\diff\tau\Bigg\|_{L^\infty_\xi}
\\
\lesssim~&  \int_0^t (\tau+1)^{-1+2p_0+(\gh+4) p_1+0.48} \Bigg[(\|\px^{-1}\vp_{j_1}\|_{H^2}+\|\partial_\xi\hat h_{j_1}\|_{L^2_\xi})\|\vp_{j_2}\|_{L^2}\|\vp_{j_3}\|_{L^{\infty}}+\|\vp_{j_1}\|_Z\|\vp_{j_2}\|_Z\|\px\vp_{j_3}\|_{L^\infty}\Bigg]\diff\tau.
\end{aligned}
\end{equation}
Notice the dispersive estimate from Lemma \ref{lem:decay}, the right-hand-side of above inequality is integrable. 
It is also summable with respect to $j_1, j_2, j_3$ under $|j_3-j_2|\leq 1$ and $|j_3-j_1|\leq 1$. By using the bootstrap assumption,
\[
\sum\limits_{\P_2}\eqref{eqn829}\lesssim \ve_1^3\lesssim
\ve_0^2.
\]

\noindent{\it Near space resonance.} Then we estimate \eqref{spRes}.
Since the integrals are symmetric in $\eta_1$ and $\eta_2$, we can assume without loss of generality that $ k_3\geq k_4$.
Using integration by parts, we write the integral as
\begin{equation*}
\begin{aligned}
& \iint_{\R^2}   \frac{\Tb_1(\eta_1, \eta_2, \xi - \eta_1 - \eta_2)}{3 i t (|\eta_1|^{2} - |\xi-\eta_1-\eta_2|^{2})}\partial_{\eta_1} e^{i  t \Phi(\xi,\eta_1,\eta_2)} \hat h_{j_1}^+(\eta_1)\hat h_{j_2}^+(\eta_2) \hat h_{j_3}(\xi-\eta_1-\eta_2)\\
&\qquad \cdot [1-\psi_{j_1-3}(\eta_1+\eta_2-2\xi)]\bigg[\psi_{k_3}\bigg(\eta_1-\frac\xi3\bigg)\psi_{\leq k_3}\bigg(\eta_2-\frac\xi3\bigg) \bigg] \diff{\eta_1} \diff{\eta_2},\\
  =~&- \frac{W_1+W_2+W_3+W_4}{3i t},
\end{aligned}
\end{equation*}
where
\begin{align*}
W_1(t, \xi)=&\iint_{\R^2}   \partial_{\eta_1}\left[\frac{\Tb_1(\eta_1, \eta_2, \xi - \eta_1 - \eta_2)}{|\eta_1|^{2} - |\xi-\eta_1-\eta_2|^{2}}\right] e^{i  t \Phi(\xi,\eta_1,\eta_2)}  \hat h_{j_1}^+(\eta_1)\hat h_{j_2}^+(\eta_2) \hat h_{j_3}(\xi-\eta_1-\eta_2) \\*
& \hspace{2in} \cdot [1-\psi_{j_1-3}(\eta_1+\eta_2-2\xi)]\psi_{k_3}\bigg(\eta_1-\frac\xi3\bigg)\psi_{\leq k_3}\bigg(\eta_2-\frac\xi3\bigg) \diff{\eta_1} \diff{\eta_2},\\[1ex]
   W_2(t, \xi)=& \iint_{\R^2}   \left[\frac{\Tb_1(\eta_1, \eta_2, \xi - \eta_1 - \eta_2)}{|\eta_1|^{2} - |\xi-\eta_1-\eta_2|^{2}}\right] e^{i  t \Phi(\xi,\eta_1,\eta_2)} \partial_{\eta_1} \hat h_{j_1}^+(\eta_1)\hat h_{j_2}^+(\eta_2) \hat h_{j_3}(\xi-\eta_1-\eta_2) \\*
   & \hspace{2in} \cdot [1-\psi_{j_1-3}(\eta_1+\eta_2-2\xi)]\psi_{k_3}\bigg(\eta_1-\frac\xi3\bigg)\psi_{\leq k_3}\bigg(\eta_2-\frac\xi3\bigg) \diff{\eta_1} \diff{\eta_2},\\[1ex]
  W_3(t, \xi)=& \iint_{\R^2}   \left[\frac{\Tb_1(\eta_1, \eta_2, \xi - \eta_1 - \eta_2)}{|\eta_1|^{2} - |\xi-\eta_1-\eta_2|^{2}}\right] e^{i  t \Phi(\xi,\eta_1,\eta_2)}  \hat h_{j_1}^+(\eta_1)\hat h_{j_2}^+(\eta_2)\partial_{\eta_1} \hat h_{j_3}(\xi-\eta_1-\eta_2) \\*
  & \hspace{2in} \cdot [1-\psi_{j_1-3}(\eta_1+\eta_2-2\xi)]\psi_{k_3}\bigg(\eta_1-\frac\xi3\bigg)\psi_{\leq k_3}\bigg(\eta_2-\frac\xi3\bigg)  \diff{\eta_1} \diff{\eta_2},\\[1ex]
  W_4(t, \xi)=& \iint_{\R^2}   \left[\frac{\Tb_1(\eta_1, \eta_2, \xi - \eta_1 - \eta_2)}{|\eta_1|^{2} - |\xi-\eta_1-\eta_2|^{2}}\right] e^{i  t \Phi(\xi,\eta_1,\eta_2)}  \hat h_{j_1}^+(\eta_1)\hat h_{j_2}^+(\eta_2) \hat h_{j_3}(\xi-\eta_1-\eta_2) \\*
  & \hspace{2in} \cdot \partial_{\eta_1}\bigg[(1-\psi_{j_1-3}(\eta_1+\eta_2-2\xi))\psi_{k_3}\bigg(\eta_1-\frac\xi3\bigg)\psi_{\leq k_3}\bigg(\eta_2-\frac\xi3\bigg)\bigg]  \diff{\eta_1} \diff{\eta_2}.
\end{align*}

Then we estimate $W_1$ to $W_4$. 
%
For $W_1$, we have the symbol estimate
\begin{align*}
&\Bigg\|m(\eta_1,\eta_2,\eta_3)\psi_{j_1}(\eta_1)\psi_{j_2}(\eta_2)\psi_{j_3}(\eta_3)\upsilon_+(\eta_1)\upsilon_+(\eta_2)\psi_{k_3}\bigg(\eta_1-\frac{\eta_1+\eta_2+\eta_3}3\bigg)\psi_{\leq k_3}\bigg(\eta_2-\frac{\eta_1+\eta_2+\eta_3}3\bigg)\Bigg\|_{S^\infty}\\
\lesssim ~& 2^{-j_1}+(1+2^{2j_1})2^{-j_1-2 k_3}
.
\end{align*}
In the above estimate, we used the assumption that $k_3\ge k_4$. Also, by the triangular inequality, we bound $|2\eta_1+\eta_2-\xi|$ as follows
$$
|2\eta_1+\eta_2-\xi|
=
\left|2(\eta_1-\frac13\xi)+(\eta_2-\frac13\xi)\right|
\ge
2\left|\eta_1-\frac{\eta_1+\eta_2+\eta_3}3\right|
-
\left|\eta_2-\frac{\eta_1+\eta_2+\eta_3}3\right|
\ge
2^{k_3}.
$$
Therefore, by considering the support of $\xi$,
\begin{align}\nonumber 
\|\cutoffxi\xi (|\xi|^{\gl} + |\xi|^{\gh}) W_1\|_{L^\infty_\xi}&\lesssim 2^{j_1}(1+t)^{\gh p_1}[2^{-j_1}+(1+2^{2j_1})2^{-j_1-2k_3}]2^{k_3}\|\hat\vp_{j_1}\|_{L^\infty_\xi}\|\hat\vp_{j_2}\|_{L^\infty_\xi}\| \vp_{j_3}\|_{L^\infty}\\\label{est_W1}
&\lesssim (1+t)^{\gh p_1+2p_0}[2^{k_3}+(1+2^{2j_1})2^{-k_3}]\|\hat\vp_{j_1}\|_{L^\infty_\xi}\|\hat\vp_{j_2}\|_{L^\infty_\xi}\| \px\vp_{j_3}\|_{L^\infty}.
\end{align}

For the terms $W_2$ to $W_4$, since their symbol satisfies
\begin{align*}
&\left\|\frac{\Tb_1(\eta_1,\eta_2,\eta_3)}{\eta_1^2-\eta_3^3}\psi_{j_1}(\eta_1)\psi_{j_2}(\eta_2)\psi_{j_3}(\eta_3)\upsilon_+(\eta_1)\upsilon_+(\eta_2)\psi_{k_3}\bigg(\eta_1-\frac{\eta_1+\eta_2+\eta_3}3\bigg)\psi_{\leq k_3}\bigg(\eta_2-\frac{\eta_1+\eta_2+\eta_3}3\bigg)\right\|_{S^\infty}\\
\lesssim~& (1+2^{2j_1})2^{-j_1-k_3}, 
\end{align*}
we have the estimate for each term
\begin{align}\label{est_W2}
\|\cutoffxi\xi (|\xi|^{\gl} + |\xi|^{\gh}) W_2\|_{L^\infty_\xi}&\lesssim (1+t)^{(\gh+2) p_1+2p_0}2^{-k_3}\|\partial_{\xi}\hat h_{j_1}\|_{L^2} \left\| \vp_{j_2}\right\|_{L^2} \|\px \vp_{j_3}\|_{L^\infty},\\\label{est_W3}
\|\cutoffxi\xi (|\xi|^{\gl} + |\xi|^{\gh}) W_3\|_{L^\infty_\xi}&\lesssim (1+t)^{(\gh+2) p_1+2p_0}2^{-k_3}\|\partial_{\xi}\hat h_{j_3}\|_{L^2} \left\| \vp_{j_2}\right\|_{L^2} \| \px\vp_{j_1}\|_{L^\infty},\\\nonumber
\|\cutoffxi\xi (|\xi|^{\gl} + |\xi|^{\gh}) W_4\|_{L^\infty_\xi}&\lesssim (1+t)^{(\gh+2) p_1}2^{-k_3}(2^{-k_3}+2^{-j_1}) 2^{k_3}\|\hat\vp_{j_1}\|_{L^\infty_\xi}\|\hat\vp_{j_2}\|_{L^\infty_\xi}\| \vp_{j_3}\|_{L^\infty}\\\label{est_W4}
&\lesssim (1+t)^{(\gh+2) p_1+2p_0}(2^{-k_3}+2^{-j_1})\|\hat\vp_{j_1}\|_{L^\infty_\xi}\|\hat\vp_{j_2}\|_{L^\infty_\xi}\| \px\vp_{j_3}\|_{L^\infty}.
\end{align}

Finally, we combine the estimates from \eqref{est_W1} to \eqref{est_W4}, take the summation over $\log_2[\varrho_1(t)] \leq k_3 \leq j_1 + 5$, and use the same argument as \eqref{6.25} and \eqref{6.26}. So we have the estimate
\begin{equation}\label{eqn720}
\begin{aligned}
& \Bigg\|\cutoffxi\xi (|\xi|^{\gl} + |\xi|^{\gh})\int_0^t  \iint_{\R^2} \Tb_1(\eta_1, \eta_2, \xi - \eta_1 - \eta_2) e^{i  t \Phi(\xi,\eta_1,\eta_2)} \hat h_{j_1}^+(\eta_1)\hat h_{j_2}^+(\eta_2) \hat h_{j_3}(\xi-\eta_1-\eta_2) \\
& \hspace{1in} \cdot[1-\psi_{j_1-3}(\eta_1+\eta_2-2\xi)] \bigg[\sum\limits_{k_3, k_4\geq \log_2[\varrho_1(\tau)]}\psi_{k_3}\Big(\eta_1-\frac\xi3\Big)\psi_{k_4}\Big(\eta_2-\frac\xi3\Big) \bigg] \diff{\eta_1} \diff{\eta_2}\diff\tau\Bigg\|_{L^\infty_\xi}\\
\lesssim~&  \int_0^t (\tau+1)^{-1+2p_0+(\gh+2) p_1+0.48} \big[(\|\vp_{j_1}\|_{H^s}+\|\partial_\xi\hat h_{j_1}\|_{L^2_\xi})\|\vp_{j_2}\|_{L^2}\|\px\vp_{j_3}\|_{L^{\infty}}+\|\vp_{j_1}\|_Z\|\vp_{j_2}\|_Z\|\px\vp_{j_3}\|_{L^{\infty}}\big]\diff\tau.
\end{aligned}\end{equation}
Notice the dispersive estimate from Lemma \ref{lem:decay}, the right-hand-side of above inequality is integrable. 
It is also summable with respect to $(j_1, j_2, j_3)$ under the assumptions $j_1\leq j_2\leq j_3$ and $|j_3-j_1|\leq 2$. By using the bootstrap assumption,
\[
\sum\limits_{\P_2}\eqref{eqn720}\lesssim \ve_1^3\lesssim
\ve_0^2.
\]

\noindent{\bf The case of $(+,-)$ or $(-,+)$.} If $\eta_1>0, \eta_2<0$, in the support of the integrand, then there is space-time resonance  $(\xi, -\xi)$. We use the unit decomposition
\[
\sum\limits_{(k_1,k_2)\in\Z^2}\psi_{k_1}(\eta_1-\xi)\psi_{k_2}(\eta_2+\xi)=1,
\]
and estimate
\begin{multline*}
\iint_{\R^2} \Tb_1(\eta_1, \eta_2, \xi - \eta_1 - \eta_2) e^{i  t \Phi(\xi,\eta_1,\eta_2)} \hat h_{j_1}^+(\eta_1)\hat h_{j_2}^-(\eta_2) \hat h_{j_3}(\xi-\eta_1-\eta_2) \\
\sum\limits_{(k_1,k_2)\in\Z^2}\psi_{k_1}(\eta_1-\xi)\psi_{k_2}(\eta_2+\xi) \diff{\eta_1} \diff{\eta_2}.
\end{multline*}
If $\eta_1<0, \eta_2>0$, in the support of the integrand, then there is space-time resonance  $(-\xi, \xi)$. We use the unit decomposition
\[
\sum\limits_{(k_1,k_2)\in\Z^2}\psi_{k_1}(\eta_1+\xi)\psi_{k_2}(\eta_2-\xi)=1,
\]
and estimate
\begin{multline*}
\iint_{\R^2} \Tb_1(\eta_1, \eta_2, \xi - \eta_1 - \eta_2) e^{i  t \Phi(\xi,\eta_1,\eta_2)} \hat h_{j_1}^-(\eta_1)\hat h_{j_2}^+(\eta_2) \hat h_{j_3}(\xi-\eta_1-\eta_2) \\
\sum\limits_{(k_1,k_2)\in\Z^2}\psi_{k_1}(\eta_1+\xi)\psi_{k_2}(\eta_2-\xi) \diff{\eta_1} \diff{\eta_2}.
\end{multline*}

We can integrate by parts and estimate each integral in the same way as above, and get the same estimate as in \eqref{eqn829}.

\subsection{Space resonance estimate}\label{spres}
In this section, we estimate the integral
\begin{align}
\label{eqn8_29}
\iint_{\R^2} \Tb_1(\eta_1, \eta_2, \xi - \eta_1 - \eta_2) e^{i t \Phi(\xi,\eta_1,\eta_2)} \hat h_{j_1}(\eta_1)\hat h_{j_2}(\eta_2) \hat h_{j_3}(\xi-\eta_1-\eta_2)\psi_{k_1}(\eta_1-\xi_1)\psi_{k_2}(\eta_2-\xi_2) \diff{\eta_1} \diff{\eta_2}
\end{align}
under the conditions that
\[
|j_1-j_3|\leq 1,\qquad 
k_1 < \log_2[\varrho_1(t)], 
\qquad k_2 < \log_2[\varrho_1(t)],
\]
where $\varrho_1(t)$ is defined in \eqref{def_rho}, then sum the result over $k_1, k_2 < \log_2(\varrho_1(t))$. Notice that when $t$ is large enough, $\psi_{k_1}(\eta_1-\xi_1)\psi_{k_2}(\eta_2-\xi_2)$ will be supported on the set where $\psi_{j_1-3}(\eta_1+\eta_2-2\xi)=1$ when $(\xi_1,\xi_2)=(\xi, \xi)$; and $\psi_{k_1}(\eta_1-\xi_1)\psi_{k_2}(\eta_2-\xi_2)$ will be supported on the set where $1-\psi_{j_1-3}(\eta_1+\eta_2-2\xi)=1$ when $(\xi_1,\xi_2)=(\frac\xi3, \frac\xi3)$. So we ignored the cut-off function $\psi_{j_1-3}(\eta_1+\eta_2-2\xi)$ in the above integral.

If  $m<\log_2 \varrho_1(t)\leq {m+1}$ for $m\in \Z$, then $-\infty < k_i \le m$ for $i=1,2$, and
\[
\sum_{k_i=-\infty}^m\psi_{k_i}(\xi)\quad \text{is supported in} \quad \left\{\xi\in \R \mid |\xi|<\frac85 \cdot 2^{m}\right\}.
\]
Thus, after summing over $(k_1,k_2)$, we only need to consider \eqref{eqn8_29} with a cutoff function in the integrand
of the form
\begin{align}
\label{defofcutoff}
\cutoff(\xi,\eta_1,\eta_2,t):=\psi\bigg(\frac{\eta_1 - \xi_1}{\varrho(t)}\bigg) \cdot \psi\bigg(\frac{\eta_2 - \xi_2}{\varrho(t)}\bigg),
\end{align}
where the function $\varrho(t)$ is defined by
\begin{equation*}
\text{$\varrho(t) = 2^m$ if $2^m < \varrho_1(t) \le 2^{m+1}$}.
\end{equation*}
In particular, from \eqref{def_rho}, we have
\begin{equation}
\frac{1}{2} (t+1)^{-0.48} \le \varrho(t) < (t+1)^{-0.48}.
\label{rhodec}
\end{equation}

The points $(\xi_1, \xi_2) \in \{ (\xi / 3, \xi / 3), (\xi, \xi), (\xi, - \xi), (-\xi, \xi)\}$ are the space and space-time resonances in Sections \ref{spres}--\ref{sptmr}. We therefore need to estimate
\begin{align}
\iint_{\R^2} \Tb_1'(\eta_1, \eta_2, \xi - \eta_1 - \eta_2) e^{i A t \Phi(\xi,\eta_1,\eta_2)} \hat h_{j_1}(\eta_1)\hat h_{j_2}(\eta_2) \hat h_{j_3}(\xi-\eta_1-\eta_2)\cutoff(\xi,\eta_1,\eta_2,t)\diff{\eta_1} \diff{\eta_2}.
\end{align}
 with the cutoff function \eqref{defofcutoff} replacing $\psi_{k_1}(\eta_1-\xi_1) \psi_{k_2}(\eta_2-\xi_2)$, in which case the integral is taken over one of the following four disjoint sets
\begin{align*}
A_1&=\bigg\{(\eta_1,\eta_2)~~\bigg{|}~~ \bigg|\eta_1-\frac{\xi}3\bigg|<\frac85 \varrho(t),\quad \bigg|\eta_2-\frac{\xi}3\bigg|<\frac85 \varrho(t)\bigg\},\\
A_2&=\bigg\{(\eta_1,\eta_2)~~\bigg{|}~~ \big|\eta_1-\xi\big|<\frac85 \varrho(t),\quad \big|\eta_2-\xi\big|<\frac85 \varrho(t)\bigg\},\\
A_3&=\bigg\{(\eta_1,\eta_2)~~\bigg{|}~~\big|(\eta_1-\xi)\big|<\frac85 \varrho(t),\quad \big|\eta_2-(-\xi)\big|<\frac85 \varrho(t) \bigg\},\\
A_4&=\bigg\{(\eta_1,\eta_2)~~\bigg{|}~~\big|\eta_1-(-\xi)\big|<\frac85 \varrho(t),\quad \big|\eta_2-\xi\big|<\frac85 \varrho(t) \bigg\}.
\end{align*}
The regions $A_1$, $A_2$, $A_3$, $A_4$ are discs centered at $(\xi/3, \xi/3)$, $(\xi, \xi)$, $(\xi, -\xi)$, and $(-\xi, \xi)$, respectively. The region $A_1$ corresponds to space resonances $\xi = \xi/3 + \xi/3 + \xi/3$, while $A_2$, $A_3$, $A_4$ correspond to space-time resonances $\xi = \xi + \xi -\xi$.

When $(\eta_1,\eta_2)\in A_1$, since
\begin{equation}
\label{T1Phi}
\frac{\Tb_1(\eta_1, \eta_2, \xi - \eta_1 - \eta_2)}{\Phi(\eta_1,\eta_2, \xi) }= 
\frac{
\frac{\alpha_2}3(\eta_1^2+\eta_2^2+\xi^2+\eta_1\eta_2-\eta_1\xi-\eta_2\xi)-1
}{3(\eta_1+\eta_2)(\xi-\eta_1)(\xi-\eta_2)},
\end{equation}
and
\begin{equation}
\frac{\Tb_1(\frac\xi3, \frac\xi3, \frac\xi3)}{\Phi(\frac\xi3,\frac\xi3, \xi) }= \frac{2\alpha_2\xi^2-9}{8\xi^3},
\end{equation}
we get
\begin{equation}
\label{T1Phi}
\frac{\Tb_1(\eta_1, \eta_2, \xi - \eta_1 - \eta_2)}{\Phi(\eta_1,\eta_2, \xi) }= \frac{2\alpha_2\xi^2-9}{8\xi^3}+
O\bigg(\left|\eta_1-\frac\xi3\right|^2+\left|\eta_2-\frac\xi3\right|^2\bigg).
\end{equation}
For $m\in \Z$, let $t_m = 2^{-{m}/{0.48}}-1$ denote the time such that $\log_2 \varrho_1(t_m)=m$, and for $t\in [0,\infty)$, let $M(t)\in \Z$
be the negative
 integer such that $M(t) < \log_2 \varrho_1(t) \le M(t)+1$. Then $\varrho(t)$ and the cut-off function $\cutoff(\xi,\eta_1,\eta_2,t)$ in \eqref{defofcutoff} are discontinuous at $t=t_m$.
 After writing
\begin{align*}
e^{i \tau \Phi(\xi, \eta_1,\eta_2)}=\frac{1}{i  \Phi(\xi, \eta_1,\eta_2)} \left[\partial_\tau e^{i\tau\Phi(\xi, \eta_1,\eta_2)}\right],
\end{align*}
and integrating by parts with respect to $\tau$ in each time interval between the time discontinuities, we obtain
\[
\begin{aligned}
&\int_1^t i\xi\iint_{\R^2}  \Tb_1(\eta_1, \eta_2, \xi - \eta_1 - \eta_2)e^{i  \tau \Phi(\xi,\eta_1,\eta_2)} \hat h_{j_1}(\eta_1, \tau)\hat h_{j_2}(\eta_2, \tau) \hat h_{j_3}(\xi-\eta_1-\eta_2, \tau) \cutoff(\xi,\eta_1,\eta_2,\tau)  \diff{\eta_1} \diff{\eta_2} \diff{\tau}\\
&\qquad =  \bigg(J_1-\int_1^t J_2(\tau)\diff{\tau}\bigg),
\end{aligned}
\]
where
\begin{align*}
& J_1 =  \iint_{\R^2}  \xi \frac{\Tb_1(\eta_1, \eta_2, \xi - \eta_1 - \eta_2)}{\Phi(\xi, \eta_1,\eta_2)} \hat h_{j_1}(\eta_1, \tau)\hat h_{j_2}(\eta_2, \tau) \hat h_{j_3}(\xi-\eta_1-\eta_2, \tau)  e^{i  \tau \Phi(\xi, \eta_1,\eta_2)}\cutoff(\xi,\eta_1,\eta_2,\tau) \diff{\eta_1} \diff{\eta_2}\Big|_{\tau=t_{M(t)}}^{\tau=t}\\
& + \sum\limits_{m=M(t)+1}^{0} \iint_{\R^2}  \xi \frac{\Tb_1(\eta_1, \eta_2, \xi - \eta_1 - \eta_2)}{\Phi(\xi, \eta_1,\eta_2)} \hat h_{j_1}(\eta_1, \tau)\hat h_{j_2}(\eta_2, \tau) \hat h_{j_3}(\xi-\eta_1-\eta_2, \tau)  e^{i  \tau \Phi(\xi, \eta_1,\eta_2)}\cutoff(\xi,\eta_1,\eta_2,\tau) \diff{\eta_1} \diff{\eta_2}\Big|_{\tau=t_m}^{\tau=t_{m - 1}},\\[1ex]
& J_2(\tau) = \xi \iint_{\R^2} \frac{\Tb_1(\eta_1, \eta_2, \xi - \eta_1 - \eta_2)}{\Phi(\xi, \eta_1,\eta_2)}  e^{i  \tau \Phi(\xi, \eta_1,\eta_2)} \partial_\tau \left[\hat h_{j_1}(\eta_1, \tau)\hat h_{j_2}(\eta_2, \tau) \hat h_{j_3}(\xi-\eta_1-\eta_2, \tau)\right]\cutoff(\xi,\eta_1,\eta_2,\tau) \diff{\eta_1} \diff{\eta_2}.
\end{align*}

For $J_1$, we have from \eqref{T1Phi} that
\[
\begin{aligned}
&\bigg|(|\xi|^{\gl}+|\xi|^{\gh})\iint_{\R^2} \cutoff(\xi,\eta_1,\eta_2,t)  \xi \frac{\Tb_1(\eta_1, \eta_2, \xi - \eta_1 - \eta_2)}{\Phi(\xi, \eta_1,\eta_2)} \hat h_{j_1}(\eta_1, \tau)\hat h_{j_2}(\eta_2, \tau) \hat h_{j_3}(\xi-\eta_1-\eta_2, \tau)   e^{i  \tau \Phi(\xi, \eta_1,\eta_2)}\diff{\eta_1} \diff{\eta_2}\bigg|\\
&\lesssim \bigg|(|\xi|^{\gl}+|\xi|^{\gh})\iint_{\R^2} \cutoff(\xi,\eta_1,\eta_2,t)  \xi^2 \hat h_{j_1}(\eta_1, \tau)\hat h_{j_2}(\eta_2, \tau) \hat h_{j_3}(\xi-\eta_1-\eta_2, \tau) e^{i\tau\Phi(\xi, \eta_1,\eta_2)} \diff{\eta_1} \diff{\eta_2}\bigg|\\
&\quad+ (|\xi|^{\gl}+|\xi|^{\gh})\iint_{\R^2} \cutoff(\xi,\eta_1,\eta_2,t)   [\varrho(\tau)]^2 \left| \hat h_{j_1}(\eta_1, \tau)\hat h_{j_2}(\eta_2, \tau) \hat h_{j_3}(\xi-\eta_1-\eta_2, \tau) \right|\diff{\eta_1} \diff{\eta_2}\\
&\lesssim \|(|\xi|^{\gl}+|\xi|^{\gh}) \hat h_{j_1}\|_{L^\infty_\xi}\|(|\xi|^{\gl}+|\xi|^{\gh}) \hat h_{j_2}\|_{L^\infty_\xi}\|(|\xi|^{\gl}+|\xi|^{\gh}) \hat h_{j_3}\|_{L^\infty_\xi}\left([\varrho(\tau)]^{2}+[\varrho(\tau)]^{4}\right).
\end{aligned}
\]
If $(\eta_1,\eta_2) \in A_1$, then the number of terms in the sum over $j_1$, $j_2$, $j_3$ is of the order $\log(t+1)$, so the right-hand side of this inequality is uniformly bounded for $\tau \geq 0$ after summing over $j_1$, $j_2$, $j_3$.

After taking the time derivative $\partial_\tau$, the term $J_2$ can be written as a sum of three terms:
\begin{align*}
&\xi \iint_{\R^2} \frac{\Tb_1(\eta_1, \eta_2, \xi - \eta_1 - \eta_2)}{\Phi(\xi, \eta_1,\eta_2)}  e^{i  \tau \Phi(\xi, \eta_1,\eta_2)}  \left[\partial_\tau\hat h_{j_1}(\eta_1, \tau)\hat h_{j_2}(\eta_2, \tau) \hat h_{j_3}(\xi-\eta_1-\eta_2, \tau)\right]\cutoff(\xi,\eta_1,\eta_2,\tau) \diff{\eta_1} \diff{\eta_2},\\[1ex]
&\xi \iint_{\R^2} \frac{\Tb_1(\eta_1, \eta_2, \xi - \eta_1 - \eta_2)}{\Phi(\xi, \eta_1,\eta_2)}  e^{i  \tau \Phi(\xi, \eta_1,\eta_2)}  \left[\hat h_{j_1}(\eta_1, \tau) \partial_\tau \hat h_{j_2}(\eta_2, \tau) \hat h_{j_3}(\xi-\eta_1-\eta_2, \tau)\right]\cutoff(\xi,\eta_1,\eta_2,\tau) \diff{\eta_1} \diff{\eta_2},
\\
&\xi \iint_{\R^2} \frac{\Tb_1(\eta_1, \eta_2, \xi - \eta_1 - \eta_2)}{\Phi(\xi, \eta_1,\eta_2)}  e^{i A\tau \Phi(\xi, \eta_1,\eta_2)}  \left[\hat h_{j_1}(\eta_1, \tau)\hat h_{j_2}(\eta_2, \tau) \partial_\tau\hat h_{j_3}(\xi-\eta_1-\eta_2, \tau)\right]\cutoff(\xi,\eta_1,\eta_2,\tau) \diff{\eta_1} \diff{\eta_2}.
\end{align*}

Using equation \eqref{h-FT-fin}, the bootstrap assumptions, and Lemma \ref{lem:decay}, we have
\begin{align*}
\|\partial_\tau\hat h\|_{L^\infty_\xi}&\lesssim \|\xi\iint_{\R^2}\Tb_1(\eta_1,\eta_2,\xi-\eta_1-\eta_2) e^{i  \tau \Phi(\xi,\eta_1,\eta_2)}\hat h(\xi-\eta_1-\eta_2)\hat h(\eta_1)\hat h(\eta_2)\diff\eta_1\diff\eta_2\|_{L^\infty_\xi}+\|\widehat{\mathcal N_{\geq 4}(\vp)}\|_{L^\infty_\xi}\\
&\lesssim\|\vp\|_{H^2}\|\px\vp\|_{W^{2,\infty}} \|\vp\|_{L^2}\\
&\lesssim \ve_1^3 (\tau + 1)^{p_0-\frac12}.
\end{align*}

Therefore, using this estimate in the $J_2$-terms and the fact that $\ve_1^3 \lesssim 
\ve_0^2
$, we get that
\[
\begin{aligned}
\left| (|\xi|^{\gl}+|\xi|^{\gh}) J_2(\tau)\right|
\lesssim \sum \| h_{\ell_1}\|_{Z}\|\partial_\tau\hat h_{\ell_2}\|_{L^\infty_\xi}\| h_{\ell_3}\|_{Z}[\varrho(\tau)]^{2}
\lesssim \ve_0 (\tau + 1)^{p_0-\frac12}[\varrho(\tau)]^{2}\sum \| h_{\ell_1}\|_{Z}\| h_{\ell_3}\|_{Z},
\end{aligned}
\]
where the summation is taken over permutations $\ell_1$, $\ell_2$, $\ell_3$ of $j_1$, $j_2$, $j_3$ for $(\eta_1,\eta_2)$ in the space-resonance region $A_1$. Again, since the number of summations is of the order
$\log(\tau + 1)$, the resulting sum is integrable over $\tau \in (1, \infty)$.

\subsection{Space-time resonances}\label{sptmr}
Considering the modified scattering and \eqref{defU}, in the following, we will estimate
\begin{align*}
&  i\xi \iint_{A_2\bigcup A_3\bigcup A_4}  \cutoff(\xi,\eta_1,\eta_2,t)\Tb_1(\eta_1, \eta_2, \xi - \eta_1 - \eta_2) e^{it\Phi(\eta_1,\eta_2, \xi)} \hat h_{j_1}(\eta_1)\hat h_{j_2}(\eta_2) \hat h_{j_3}(\xi-\eta_1-\eta_2)  \diff{\eta_1} \diff{\eta_2}\\
-&
\frac{i(2\alpha_2\xi^2-1)\pi}{6}
\frac{
|\hat{h}(t, \xi)|^2
\hat{h}(t, \xi)
}
{t}.
\end{align*}
The estimates for $A_2$, $A_3$, and $A_4$ are similar, so we only present the details for the $A_2$ integral.
The corresponding integral for $A_2$ in \eqref{modscat} can be decomposed into
\begin{align}
\label{modscat1}
\begin{split}
&i\xi \iint_{\R^2} e^{i t \Phi(\xi, \eta_1, \eta_2)} \psi\bigg(\frac{\eta_1 - \xi}{ \varrho(t)}\bigg) \cdot \psi\bigg(\frac{\eta_2 - \xi}{ \varrho(t)}\bigg)\\*
&\qquad \cdot \bigg[\Tb_1(\eta_1, \eta_2, \xi - \eta_1 - \eta_2) \hat{h}_{j_1}(\eta_1)\hat{h}_{j_2}(\eta_2) \hat{h}_{j_3}(\xi-\eta_1-\eta_2) - \Tb_1(\xi, \xi, -\xi)|\hat{h}(\xi)|^2\hat{h}(\xi)\bigg] \diff{\eta_1} \diff{\eta_2}\\
\end{split}
\end{align}
and
\begin{align}
\label{modscat2}
i\xi \Tb_1(\xi, \xi, - \xi) |\hat{h}(t,\xi)|^2\hat{h}(t,\xi) \bigg[\iint_{A_2} e^{i t \Phi(\xi, \eta_1, \eta_2)} \psi\bigg(\frac{\eta_1 - \xi}{ \varrho(t)}\bigg) \cdot \psi\bigg(\frac{\eta_2 - \xi}{ \varrho(t)}\bigg) \diff{\eta_1} \diff{\eta_2} - \frac{2 \pi }{6\xi t}\bigg].
\end{align}

The estimates for \eqref{modscat1} are achieved by a Taylor expansion,
\begin{align*}
&\bigg|(|\xi|^{\gl}+|\xi|^{\gh})i\xi \iint_{\R^2} e^{i t \Phi(\xi, \eta_1, \eta_2)} \psi\bigg(\frac{\eta_1 - \xi}{ \varrho(t)}\bigg) \cdot \psi\bigg(\frac{\eta_2 - \xi}{ \varrho(t)}\bigg)
\\*
&\qquad \cdot \Big[\Tb_1(\eta_1, \eta_2, \xi - \eta_1 - \eta_2) \hat{h}_{j_1}(\eta_1)\hat{h}_{j_2}(\eta_2) \hat{h}_{j_3}(\xi-\eta_1-\eta_2) - \Tb_1(\xi, \xi, -\xi)|\hat{h}(\xi)|^2\hat{h}(\xi) \Big] \diff{\eta_1} \diff{\eta_2}\bigg|
\\
\lesssim~& (|\xi|^{\gl}+|\xi|^{\gh}) |\xi| \iint_{\R^2} 
\psi\bigg(\frac{\eta_1 - \xi}{ \varrho(t)}\bigg) \cdot \psi\bigg(\frac{\eta_2 - \xi}{ \varrho(t)}\bigg)
\\
\Bigg[
&
\left|\partial_{\eta_1}\left[\Tb_1(\eta_1, \eta_2, \xi - \eta_1 - \eta_2) \hat{h}_{j_1}(\eta_1)\hat{h}_{j_2}(\eta_2) \hat{h}_{j_3}(\xi-\eta_1-\eta_2)\right] \bigg|_{\eta_1 = \eta_1'} (\xi-\eta_1)\right|
\\
+&\left|\partial_{\eta_2}\left[\Tb_1(\eta_1, \eta_2, \xi - \eta_1 - \eta_2) \hat{h}_{j_1}(\eta_1)\hat{h}_{j_2}(\eta_2) \hat{h}_{j_3}(\xi-\eta_1-\eta_2)\right] \bigg|_{\eta_2 = \eta_2'} (\xi-\eta_2)\right| 
\Bigg]\diff{\eta_1} \diff{\eta_2}
\end{align*}
Since 
$$
{\bf T_1}(\eta_1,\eta_2,\xi-\eta_1-\eta_2)
=
\frac{\alpha_2}3(\eta_1^2+\eta_2^2+\xi^2+\eta_1\eta_2-\eta_1\xi-\eta_2\xi)-1,
$$
$$
\partial_{\eta_1}{\bf T_1}(\eta_1,\eta_2,\xi-\eta_1-\eta_2)
=
2\eta_1+\eta_2-\xi
$$
and
$$
\partial_{\eta_2}{\bf T_1}(\eta_1,\eta_2,\xi-\eta_1-\eta_2)
=
2\eta_2+\eta_1-\xi
$$
we get 
\begin{align*}
\eqref{modscat1}
\lesssim &\|(1+|\xi|^{\frac52})\hat\vp_{j_1}\|_{L^\infty_\xi} \|(1+|\xi|^{\frac52})\hat\vp_{j_2}\|_{L^\infty_\xi} \|(1+|\xi|^{\frac52})\hat\vp_{j_3}\|_{L^\infty_\xi}[\varrho(t)]^{3}\\
&+ \sum\|(|\xi|^{\frac12}+|\xi|^{\frac52})\hat\vp_{\ell_1}\|_{L^\infty_\xi}\|(|\xi|^{\frac12}+|\xi|^{\frac52})\hat\vp_{\ell_2}\|_{L^\infty_\xi}\|\S\vp_{\ell_3}\|_{L^2}[\varrho(t)]^{5/2}(t+1)^{2p_1},
\end{align*}
where $\ell_1$, $\ell_2$ and $\ell_3$ are the permutation of $j_1,j_2$ and $j_3$. The sum is with respect to the permutation.
As for \eqref{modscat2}, it suffices to estimate
\begin{align*}
& \bigg|(|\xi|^{\gl}+|\xi|^{\gh}) i\xi \Tb_1(\xi, \xi, - \xi) |\hat{h}(t,\xi)|^2\hat{h}(t,\xi) \bigg[\iint_{\R^2} e^{i t \Phi(\xi, \eta_1, \eta_2)} \psi\bigg(\frac{\eta_1 - \xi}{ \varrho(t)}\bigg) \cdot \psi\bigg(\frac{\eta_2 - \xi}{ \varrho(t)}\bigg) \diff{\eta_1} \diff{\eta_2} - \frac{2 \pi }{6\xi t}\bigg]\bigg|
\\
\lesssim ~& \|(|\xi|^{\gl}+|\xi|^{\gh}) \hat{\vp}(\xi)\|_{L^\infty_\xi} \|\hat{\vp}(\xi)\|_{L^\infty_\xi} \|(1+|\xi|^2)\hat{\vp}(\xi)\|_{L^\infty_\xi}
\\
& \qquad \cdot \bigg\|
\xi
\iint_{\R^2} e^{i t \Phi(\xi, \eta_1, \eta_2)} \psi\bigg(\frac{\eta_1 - \xi}{ \varrho(t)}\bigg) \cdot \psi\bigg(\frac{\eta_2 - \xi}{ \varrho(t)}\bigg) \diff{\eta_1} \diff{\eta_2} - \frac{2 \pi }{6 t}\bigg\|_{L^\infty_\xi}
\\
\lesssim ~& \|\vp\|_{Z}^3 \bigg\|
\xi
\iint_{\R^2} e^{i t \Phi(\xi, \eta_1, \eta_2)} \psi\bigg(\frac{\eta_1 - \xi}{ \varrho(t)}\bigg) \cdot \psi\bigg(\frac{\eta_2 - \xi}{ \varrho(t)}\bigg) \diff{\eta_1} \diff{\eta_2} - \frac{2 \pi }{6 t}\bigg\|_{L^\infty_\xi}.
\end{align*}
Writing $(\eta_1, \eta_2) = (\xi + \zeta_1, \xi + \zeta_2)$, we find from \eqref{defPhi} that
\begin{align*}
\Phi(\xi, \eta_1, \eta_2)
=
3(\eta_1+\eta_2)(\xi-\eta_1)(\xi-\eta_2)
=
6\xi\zeta_1 \zeta_2 + O\bigg(\zeta_1^3 + \zeta_2^3\bigg)
=
6\xi
\zeta_1 \zeta_2  + O\left([ \varrho(t)]^3 \right).
\end{align*}
Since  $ \varrho(t)$ satisfies \eqref{rhodec}, the error term is integrable in time, so
we now only need to estimate
\begin{align}
\label{modscat}
\begin{split}
J_3 &:= \bigg\|
\xi
\iint_{\R^2} e^{- i t 6\xi \zeta_1 \zeta_2} \psi\bigg(\frac{\zeta_1}{ \varrho(t)}\bigg) \cdot \psi\bigg(\frac{\zeta_2}{ \varrho(t)}\bigg) \diff{\zeta_1} \diff{\zeta_2} - \frac{2 \pi }{6 t}\bigg\|_{L^\infty_\xi}.
\end{split}
\end{align}
Making the change of variables
\[
\zeta_1 = \sqrt{\frac{1}{6\xi t}} x_1,\qquad \zeta_2 = \sqrt{\frac{1}{6\xi t}} x_2,
\]
in \eqref{modscat}, we obtain
\begin{align}
\label{modscatest}
\begin{split}
J_3 &\le  
\frac{1}{6t}
\bigg\|\iint_{\R^2} e^{- i x_1 x_2} \psi\bigg(\frac{1}{\sqrt{6\xi t}  \varrho(t)} x_1\bigg) \cdot \psi\bigg(\frac{1}{\sqrt{6\xi t}  \varrho(t)} x_2\bigg) \diff{x_1} \diff{x_2} - 2 \pi\bigg\|_{L^\infty_\xi}.
\end{split}
\end{align}
The integral identity
\[
\int_\R e^{- a x^2 - b x} \diff{x} =\sqrt{ \frac{\pi}{a}} e^{\frac{b^2}{4 a}}, \qquad \text{for all $a, b \in \mathbb{C}$ with $\Re{a} > 0$}
\]
gives that
\[
\iint_{\R^2} e^{- i x_1 x_2} e^{- \frac{x_1^2}{\mathfrak{B}^2}} e^{- \frac{x_2^2}{\mathfrak{B}^2}} \diff{x_1} \diff{x_2} = \sqrt{\pi} \mathfrak{B} \int_\R e^{- \frac{x_2^2}{\mathfrak{B}^2}} e^{- \frac{\mathfrak{B}^2 x_2^2}{4}} \diff{x_2} = 2 \pi + O(\mathfrak{B}^{-1}),
\qquad \text{as $\mathfrak{B}\to \infty$},
\]
and therefore
\begin{align}
\label{2piint}
\left|\iint_{\R^2} e^{- i x_1 x_2} \psi\bigg(\frac{x_1}{\mathfrak{B}}\bigg) \psi\bigg(\frac{x_2}{\mathfrak{B}}\bigg) \diff{x_1} \diff{x_2} -2 \pi\right|\lesssim \mathfrak{B}^{- 1 / 2},
\qquad \text{as $\mathfrak{B}\to \infty$}.
\end{align}
Using \eqref{2piint} with $\mathfrak{B} = \sqrt{6\xi t}  \varrho(t)  = O(t^{0.02})$ in \eqref{modscatest}, then yields
\begin{align*}
J_3 \lesssim (t + 1)^{-1.01}.
\end{align*}
Hence, the right-hand side decays faster in time than $1 / t$, which implies that \eqref{modscat2} is integrable in time and bounded by a constant multiple of $\ve_0^3$. 

Putting all the above estimates together, we conclude that
\[
\int_0^{\infty} \|(|\xi|^{\gl}+|\xi|^{\gh}) U(t, \xi)\|_{L^\infty_\xi} \diff t\lesssim \ve_0.
\]

\subsection{Higher-degree terms}\label{higherdegree}
In this subsection, we prove that
\[
\left\|\cutoffxi(|\xi|^{\gl}+|\xi|^{\gh})\widehat{\mathcal{N}_{\geq4}(\vp)}\right\|_{L^\infty_\xi}
\]
is integrable in time.  

Using the support of $\cutoffxi$, we have
\begin{align*}
\left\|\cutoffxi(|\xi|^{\gl}+|\xi|^{\gh})\widehat{\mathcal{N}_{\geq4}(\vp)}\right\|_{L^\infty_\xi}
\lesssim&(t+1)^{4p_1}\|\px^{-1}\Nc_{\geq4}\|_{L^1}\\
\lesssim&(t+1)^{4p_1}(\|\px^{-1}\Nc_{4}\|_{L^1}+\|\px^{-1}\Nc_{\geq5}\|_{L^1}).
\end{align*}
For $\Nc_{\geq5}(\vp)$, by using bootstrap assumption, we have
\begin{align*}
\|\px^{-1}\Nc_{\geq5}\|_{L^1}\lesssim \||\vp|_2|\px\vp|_2\|_{L^\infty}\|\vp\|_{L^{\infty}}^2\|\vp\|_{L^2}\lesssim (t+1)^{-1-\frac23}\ve_1^5,
\end{align*}
which is integrable with $t\in(0,\infty)$. 

For $\Nc_4(\vp)$, 
\begin{align*}
\F[\Nc_4(\vp)]=&i\xi\iiint \Tb_2(\eta_1, \eta_2, \eta_3, \xi-\eta_1-\eta_2-\eta_3)\hat\vp(\eta_1)\hat\vp(\eta_2)\hat\vp(\eta_3)\hat\vp(\xi-\eta_1-\eta_2-\eta_3)d\eta_1d\eta_2d\eta_3\\
=&i\xi\iiint \Tb_2(\eta_1, \eta_2, \eta_3, \xi-\eta_1-\eta_2-\eta_3)e^{i\Psi t}\hat h(\eta_1)\hat h(\eta_2)\hat h(\eta_3)\hat h(\xi-\eta_1-\eta_2-\eta_3)d\eta_1d\eta_2d\eta_3
\end{align*}
where the symbol 
\[
\Tb_2(\eta_1, \eta_2, \eta_3,\eta_4)=-2\eta_1^2-2\eta_1\eta_2-\eta_1\eta_3,
\]
and the phase
\[
\Psi=\Psi(\eta_1,\eta_2,\eta_3,\xi-\eta_1-\eta_2-\eta_3)=\xi^3-(\xi-\eta_1-\eta_2-\eta_3)^3-\eta_1^3-\eta_2^3-\eta_3^3.
\]
Notice that 
\[
\partial_{\eta_i}\Psi=3(\xi-\eta_1-\eta_2-\eta_3)^2-3\eta_i^2, \quad i=1,2,3
\]
$|\nabla_{\etab}\Psi|=0$ implies $|\eta_1|=|\eta_2|=|\eta_3|=|\xi-\eta_1-\eta_2-\eta_3|$. The space-time resonance corresponds to $\Psi=|\nabla_{\etab}\Psi|=0$, which does not have non-zero solutions. So, there is no space-time resonance. By similar estimates as of Sections \ref{nonre}-\ref{spres}, we have 
\[
\int_0^\infty\left\|\cutoffxi(|\xi|^{\gl}+|\xi|^{\gh})\widehat{\mathcal{N}_{4}(\vp)}\right\|_{L^\infty_\xi} dt\lesssim\ve_1^4.
\]

This completes the proof of Proposition~\ref{Bootstrap-prop}.

\appendix

\end{document}